\documentclass{article}

\usepackage{arxiv}

\usepackage[utf8]{inputenc} 
\usepackage[T1]{fontenc}    
\usepackage{hyperref}       
\hypersetup{
	colorlinks=true,
	citecolor=blue
}
\usepackage{url}            
\usepackage{booktabs,multirow}       
\usepackage{amsfonts,amsmath,amssymb,amsthm}       
\usepackage{nicefrac}       
\usepackage{microtype}      
\usepackage{cleveref}       
\usepackage{graphicx}
\graphicspath{{../figures/}}
\usepackage{natbib}
\usepackage{enumitem}

\numberwithin{equation}{section}
\theoremstyle{plain}

\newtheorem{theorem}{Theorem}[section]

\newtheorem{corollary}{Corollary}

\theoremstyle{remark}
\newtheorem{definition}[theorem]{Definition}

\newtheorem{remark}{Remark}[section]

\def\E{\mathbb{E}}
\def\R{\mathbb{R}}

\newcommand{\bb}[1]{\boldsymbol{#1}}
\newcommand{\bbhat}[1]{\widehat{\boldsymbol{#1}}}
\newcommand{\tr}{^{\intercal}}

\newcommand{\norm}[1]{\left\Vert #1 \right\Vert}   

\DeclareMathOperator*{\argmin}{arg\,min}
\newcommand{\dash}{^{\prime}}

\newcommand{\cov}{\mathrm{Cov}}
\newcommand{\var}{\mathrm{Var}}
\newcommand{\normdist}{\mathcal{N}}

\newcommand{\Hcal}{\mathcal{H}}
\newcommand{\Fcal}{\mathcal{F}}
\newcommand{\Hdiff}[1]{\partial_{#1}\mathcal{H}_{n,p}}
\newcommand{\Hdifftwo}[1]{\partial^2_{#1}\mathcal{H}_{n,p}}

\title{Analysis of the rSVDdpd Algorithm: A Robust Singular Value Decomposition Method using Density Power Divergence}

\author{Subhrajyoty Roy\\
	Interdisciplinary Statistical Research Unit\\
	Indian Statistical Insititute, Kolkata \\
	\texttt{roysubhra98@gmail.com} \\
	\And
	Abhik Ghosh\\
	Interdisciplinary Statistical Research Unit\\
	Indian Statistical Insititute, Kolkata \\
	\texttt{abhik.ghosh@isical.ac.in} \\
	\And
	Ayanendranath Basu\\
	Interdisciplinary Statistical Research Unit\\
	Indian Statistical Insititute, Kolkata \\
	\texttt{ayanbasu@isical.ac.in} \\
}

\renewcommand{\shorttitle}{Theoretical Analysis of rSVDdpd}

\hypersetup{
pdftitle={Analysis of the rSVDdpd Algorithm: A Robust Singular Value Decomposition Method using Density Power Divergence},
pdfauthor={Subhrajyoty Roy, Abhik Ghosh, Ayanendranath Basu},
pdfkeywords={Singular value decomposition, Matrix factorization, Density power divergence},
}

\begin{document}
\maketitle

\begin{abstract}
    The traditional method of computing singular value decomposition (SVD) of a data matrix is based on a least squares principle, thus, is very sensitive to the presence of outliers. Hence the resulting inferences across different applications using the classical SVD are extremely degraded in the presence of data contamination (e.g., video surveillance background modelling tasks, etc.). A robust singular value decomposition method using the minimum density power divergence estimator (rSVDdpd) has been found to provide a satisfactory solution to this problem and works well in applications. For example, it provides a neat solution to the background modelling problem of video surveillance data in the presence of camera tampering. In this paper, we investigate the theoretical properties of the rSVDdpd estimator such as convergence, equivariance and consistency under reasonable assumptions. Since the dimension of the parameters, i.e., the number of singular values and the dimension of singular vectors can grow linearly with the size of the data, the usual M-estimation theory has to be suitably modified with concentration bounds to establish the asymptotic properties. We believe that we have been able to accomplish this satisfactorily in the present work. We also demonstrate the efficiency of rSVDdpd through extensive simulations.
\end{abstract}


\section{Introduction}\label{sec:intro}

Singular value decomposition (SVD) is a matrix factorization method that breaks down a real or complex matrix into three parts: two orthogonal matrices consisting of the singular vectors and one diagonal matrix with non-negative diagonal entries made up of singular values. SVD is commonly viewed as a low-rank approximation of a linear transformation, representing the transformation as a sequence of rotations and dilations. The applications of SVD are diverse. Its mathematical applications include computing pseudoinverses or Moore-Penrose inverses of matrices, efficiently solving systems of homogeneous linear equations, determining range space, null space, and rank of a linear transformation, and finding ordinary least square solutions using the Golub and Reinsch algorithm~\citep{golub1971singular}. SVD has also been widely utilized in statistical and machine-learning methods for data analysis and modelling in various real-world applications. One notable application is principal component analysis (PCA), which employs SVD to decompose a data matrix into a lower-dimensional representation that captures the maximum variability in the original data using a reduced set of variables. Other popular dimension reduction techniques  such as correspondence analysis~\citep{greenacre2017correspondence}, latent semantic indexing~\citep{hofmann1999probabilistic,anandarajan2019semantic} and clustering techniques~\citep{drineas2004clustering,cheng2019bearing} also rely on SVD as a fundamental component. SVD is extensively employed in pattern recognition within signal, image, and video processing domains. Its applications include image watermarking schemes~\citep{dappuri2020non}, signal denoising and feature enhancements~\citep{zhao2017novel}, audio watermarking~\citep{rezaei2019robust}, sound source localization~\citep{grondin2019multiple}, and sound recovery techniques~\citep{zhang2016note}. Moreover, SVD has gained prominence in the field of bioinformatics, where it is used for analyzing protein functional associations~\citep{franceschini2016svd}, clustering gene expression data~\citep{bustamam2018clustering}, and predicting protein-coding regions~\citep{das2017proteincoding}. In the realm of geographical science, \cite{kumar2011new} employed SVD-based techniques, including its robust variant, to generate accurate graphical representations of climate data, mitigating the impact of extreme weather events such as thunderstorms and heavy rainfall. Such a wide range of applications clearly underscore the relevance of SVD as an extremely integral component of data analysis across a multitude of disciplines.

However, as already indicated by several authors~\citep{niss2001};~\citep{liu2003robust};~\citep{kumar2011new}, the usual method of computing SVD is highly susceptible to outliers present in the data matrix. In contrast, as the data are becoming increasingly vast and complex in the recent era, it is also susceptible to the inclusion of different forms of noises, corruptions and contamination by outlying observations. This readily introduces the need for a robust algorithm for the computation of SVD which remains unaffected (or is minimally affected) due to the presence of outliers, and thus leads to more stable and trustworthy solutions in different applications mentioned previously. 

\subsection{Related Works and Our Contribution}\label{sec:related-works}

\cite{ammann1993} was one of the early pioneers in developing a robust version of the SVD. He treated it as a special case of the projection pursuit problem to be solved using the transposed QR algorithm. Other researchers, such as~\cite{niss2001,liu2003robust} and \cite{qifa2005}, approached the computation of SVD as a least squares problem and proposed robust extensions using alternating $L_1$ regression algorithms with the least absolute deviation (LAD) loss function. However, a simple LAD approach is sensitive to high leverage points, which led to the exploration of weighted LAD approaches and the utilization of the Huber weight function~\citep{jung-wlad}. In a different attempt, \cite{rey2007-total-svd} introduced a robust method called ``Total'' SVD, which employed Huber's weight function and ``Total'' least squares~\citep{markovsky2007overview}. This approach accounted for errors in both the data matrix and the singular vectors, in contrast to the usual least squares where the only source of error is the response variable. Although this resulted in a more robust SVD estimate, the method faced several convergence issues as mentioned by~\cite{rey2007-total-svd} himself. Alternatively, \cite{zhang2013} incorporated the Huber weight function in the loss function and combined it with a squared error-based penalty function for regularization, creating another robust SVD estimator. \cite{Wang_2017} used an estimator derived from an $\alpha$-stable distribution with a cost function $\rho(x) = \log(x^2 + K^2)$, where the tuning parameter $K$ provides a balance between robustness and efficiency. Nevertheless, finding the appropriate tuning parameter $K$ for the estimator was challenging. Apart from a simple alternating $L_1$ regression approach~\citep{gabriel1979lower}, there is a lack of theoretical guarantees or properties for the resulting SVD estimates in the literature. Convergence and orthogonality of the singular vectors are not assured, and the computational complexity poses a significant challenge, especially for large data matrices.

\cite{roy2021new} introduced the rSVDdpd (Robust Singular Value Decomposition using Density Power Divergence) estimator as an outlier-resistant matrix factorization technique in the context of video surveillance background modelling. They transformed the SVD computation into an alternating regression problem and utilized the density power divergence (DPD) loss function~\citep{basu-dpd}. The minimum DPD estimator (MDPDE) has exhibited strong robustness and high efficiency in statistics and information theory (see, e.g., \cite{basu-dpd}; a brief description is also provided in Appendix~\ref{appendix:dpd}). The use of the DPD loss function in the proposed rSVDdpd estimator benefits from these desirable statistical properties. This paper aims to provide theoretical justifications for the rSVDdpd algorithm by establishing its convergence and mathematical properties like equivariance and asymptotic consistency. The primary challenge here is to extend the existing results on the MDPDE type estimators~\citep{ghosh2013} to the case where the dimension of the parameter grows to infinity linearly in sample size, and thus, concentration bounds are needed to ensure that the desirable asymptotic properties hold. We also conduct extensive simulation studies to demonstrate its applicability as a general-purpose robust matrix factorization technique and its performance compared to the existing approaches by~\cite{zhang2013} and~\cite{niss2001}, where rSVDdpd outperforms them in most simulation setups.

\section{Description of the rSVDdpd estimator}\label{sec:description}

The problem of Singular Value Decomposition starts with a data matrix $\bb{X}$ of dimension $n \times p$ ($n$ and $p$ may be different), admitting an approximate low-rank representation of the form
\begin{equation}
    \bb{X} = \sum_{k = 1}^{r} \lambda_k \bb{u}_{k}\bb{v}_{k}\tr + \bb{E},
    \label{eqn:rsvd-decomp}
\end{equation}
\noindent where $\{\bb{u}_{k}\}_{k=1}^r$ is a set of $r$ orthonormal vectors of length $n$ and $\{\bb{v}_{k}\}_{k=1}^r$ is a set of $r$ orthonormal vectors of length $p$. The $n \times p$ dimensional matrix $\bb{E}$ consists of the errors $e_{ij}$s, which are generally expected to be smaller in magnitude than the corresponding entries of the data matrix $\bb{X}$, except at a few coordinates with outlying observations. The goal is to estimate the unknown rank $r$ of the low-rank component of $\bb{X}$, the left and right singular vectors $\bb{u}_k$s and $\bb{v}_k$s and the nonnegative singular values $\lambda_k$s. For notational convenience, let us denote $\bb{U} = [\bb{u}_1, \dots \bb{u}_r]_{n \times r}$, $\bb{V} = [\bb{v}_1, \dots \bb{v}_r]_{p \times r}$ and $\bb{\Lambda}$ as the $r\times r$ diagonal matrix with diagonal entries $\lambda_1, \lambda_2, \dots, \lambda_r$. For the time being, we assume that the rank $r$ is known and focus on estimating the singular values and vectors.

The above description of SVD as in~\eqref{eqn:rsvd-decomp} is equivalent to the LSN decomposition in~\cite{zhou2010stable} and is a generalization of the LS decomposition of~\cite{candes2011}. Because of the presence of outlying values in the errors $e_{ij}$s, one can consider them as independent random variables with $e_{ij}$ following a mixture distribution of the form $G_{ij} = (1-\delta)G_{1,ij} + \delta G_{2,ij}$ for some small $\delta \in [0, 1]$, denoting the proportion of contamination. Here, $G_{1,ij}$ and $G_{2,ij}$ are the distribution functions corresponding to the dense perturbation and the sparse outlying components of $e_{ij}$. We assume that the distribution $G_{ij}$ admits a density $g_{ij}$ with respect to the Lebesgue measure for all $i = 1, \dots , n$ and $j = 1, \dots, p$. On the flip side, if the estimated singular values and the vectors are correctly estimated, then the errors may be modelled as independent and identically distributed (i.i.d.) observations from a symmetrically distributed scale family of densities, $\Fcal = \{ \sigma^{-1} f(\cdot/\sigma): \sigma \in (0, \infty) \}$ with $f(x) = f(-x)$ for all $x \in \R$, where the functional form of $f$ is known. A popular and standard choice of $f$ may be the standard normal density function. Hence, the problem of estimating SVD robustly as in decomposition~\eqref{eqn:rsvd-decomp} can be regarded as a robust estimation problem. To solve this, \cite{roy2021new} use the minimum density power divergence estimator (MDPDE)~\citep{basu-dpd} which have been shown to possess strong robustness properties. In this case of low-rank decomposition as in~\eqref{eqn:rsvd-decomp}, the MDPDE is given by the minimizer of
\begin{equation}
    H_\alpha^{(r)}(\bb{\theta}) = \dfrac{1}{np}\sum_{i=1}^n \sum_{j=1}^p V_{ij,\alpha}^{(r)}(\bb{\theta})
    \label{eqn:H-general}
\end{equation}
\noindent where
\begin{equation}
    V_{ij,\alpha}^{(r)}(\bb{\theta}) = \sigma^{-\alpha} \left[ M_f - \left(1+\dfrac{1}{\alpha}\right)\right.
        \left. f^\alpha\left( \left\vert \dfrac{X_{ij} - \sum_{k=1}^r \lambda_k u_{ki}v_{kj}}{\sigma} \right\vert \right) \right]
    \label{eqn:V-general}
\end{equation}
\noindent and $M_f = \int f^{1+\alpha}$. Here, the parameter $\bb{\theta} = \left( \bb{\Lambda}, \bb{U}, \bb{V}, \sigma^2 \right)$ is restricted in the parameter space $[0, \infty)^r \times S_n^r \times S_p^r \times (0,\infty)$, where $S_n^r$ and $S_p^r$ denote the $r$-Stiefel manifolds of order $n$ and $p$ respectively. The resulting matrices $\widehat{\bb{\Lambda}}$, $\widehat{\bb{U}}$ and $\widehat{\bb{V}}$ as a solution to the MDPDE objective function given in~\eqref{eqn:H-general} is then defined to be the Robust SVD using Density Power Divergence (rSVDdpd) estimator of the data matrix $\bb{X}$ up to rank $r$.

A standard and popular choice for the scale family of densities is to consider the normal densities with mean $0$ and unknown variance $\sigma^2$. In this case, the $V$-function as in~\eqref{eqn:V-general} reduces to
\begin{equation}
    V_{ij,\alpha}^{(r)}(\bb{\theta}) = \sigma^{-\alpha}\left[ \dfrac{1}{\sqrt{1+\alpha}}  \right.\\
        \left. - \left(1 + \dfrac{1}{\alpha} \right) e^{-\alpha (X_{ij} - \sum_{k=1}^r \lambda_k u_{ki}v_{kj})^2 / 2\sigma^2} \right].
    \label{eqn:V-normal}
\end{equation}

A direct minimization of the objective function given in~\eqref{eqn:H-general} is extremely difficult to solve since the quantities $\bb{U}$ and $\bb{V}$ are restricted to Stiefel manifolds which are nonlinear and nonconvex spaces. Following the footsteps of~\cite{rey2007-total-svd}, we reformulate the decomposition in~\eqref{eqn:rsvd-decomp} as
\begin{equation}
    \bb{X} = \sum_{k = 1}^{r} \bb{a}_{k}\bb{b}_{k}\tr + E,
    \label{eqn:svd}
\end{equation}
\noindent where $\bb{a}_k$s and $\bb{b}_k$s are still orthogonal sets of vectors for $k = 1, 2, \dots r$, but not necessarily normalized. Once the estimates of $\bb{a}_k$s and $\bb{b}_k$s are known, they can be normalized to obtain the $\bb{u}_k$s and $\bb{v}_k$s and the singular values are then given by $\lambda_k = \Vert \bb{a}_k \Vert \Vert \bb{b}_k \Vert$ for each $k=1, \ldots, r$, where $\Vert\cdot\Vert$ denotes the usual Euclidean ($L_2$) norm. Equipped with this idea, \cite{roy2021new} describe the rSVDdpd estimator for a rank one decomposition of the form
\begin{equation}
    X_{ij} = a_i b_j + e_{ij}, \qquad i = 1, 2, \dots n.
    \label{eqn:regression-reduced}
\end{equation}
\noindent Following~\cite{rey2007-total-svd}, they view~\eqref{eqn:regression-reduced} as a linear regression equation instead of a singular value decomposition problem. For a fixed index $j$,~\eqref{eqn:regression-reduced} can be interpreted as a linear regression with $X_{ij}$s as an observed response, $a_i$s as the covariate values, $b_j$s as the regression slope parameters to be estimated and $e_{ij}$s as the random error component. As we vary the column index $j = 1, 2, \dots p$, we are posed with $p$ such linear regression problems, solving each of them will jointly yield an estimate of $\bb{b}=(b_1, \ldots, b_p)$ given the values of $a_i$s. Now, one can interchange the role of $a_i$s and $b_j$s and view~\eqref{eqn:regression-reduced} as regression of $X_{ij}$ on $b_j$s for fixed $i$, and estimate $\bb{a}=(a_1, \ldots, a_n)$. Then, one can alternatively estimate $\bb{a}$ and $\bb{b}$ by solving these linear regression problems, and the converged values of $(\bb{a}, \bb{b})$ will yield the required decomposition. \cite{roy2021new} aim to solve these regression problems using the popular minimum density power divergence estimator (MDPDE) introduced in~\cite{ghosh2013}. 

Denoting $\psi(x) = -f^\alpha(\vert x\vert)u(\vert x\vert)/\vert x\vert$ where $u(\cdot)$ is the score function of the density $f$, i.e., $u(x) := f'(x)/f(x)$, assuming $f'$ exists, the solutions to the alternating regression problems give rise to the iterative equations
\begin{align}
    a_i      & = \dfrac{\sum_j b_j X_{ij} \psi( e_{ij}/\sigma ) }{ \sum_j b_j^2 \psi( e_{ij} / \sigma) }, \qquad i = 1, \ldots n,
    \label{eqn:general-1} \\
    b_j      & = \dfrac{\sum_i a_i X_{ij} \psi( e_{ij}/\sigma ) }{ \sum_i a_i^2 \psi( e_{ij}/\sigma ) }, \qquad j = 1, \ldots p,
    \label{eqn:general-2} \\
    \sigma^2 & = \dfrac{(np)^{-1}\sum_i \sum_j e_{ij}^2 \psi(e_{ij} / \sigma )}{ (np)^{-1} \sum_i \sum_j \psi(e_{ij}/\sigma) - \frac{\alpha}{1+\alpha} M_f }.
    \label{eqn:general-3}
\end{align}

\begin{remark}
    If $f$ is the standard normal density, then $\psi(x) = e^{-\alpha x^2/2}$ for all $x > 0$, which leads to the exact iteration steps mentioned in~\cite{roy2021new}. For $\alpha > 0$, this is a decreasing function and it leads to a robust SVD estimator. However, for $\alpha = 0$, $\psi(x) = 1$ and the iterative equations~\eqref{eqn:general-1}-\eqref{eqn:general-3} reduces to the estimation procedure of classical SVD~\citep{niss2001}. Hence, the proposed algorithm produces a class of SVD estimators including both robust and non-robust estimators.
\end{remark}

Given that such rank-one decompositions as in~\eqref{eqn:regression-reduced} can be obtained, one can stack the estimates of $a_i$s to get $\widehat{\bb{a}}_1$, and combine estimates of $b_j$s to get $\widehat{\bb{b}}_1$, the unnormalized first set of the singular vectors. Then these can be normalized and one obtains an estimate of the first singular value $\widehat{\lambda}_1 = \Vert \widehat{\bb{a}}_1 \Vert \Vert \widehat{\bb{b}}_1\Vert$. For subsequent singular values, one can apply the same estimation algorithm on the residual matrix $\bb{X} - \widehat{\lambda}_1\widehat{\bb{a}}_1\widehat{\bb{b}}_1\tr$. Such an iterative method of rank one approximation is quite common in the matrix factorization literature~\citep{niss2001,chicoki-robust-nmf}. However, this method does not guarantee that the subsequent singular vectors remain orthonormal to the previous set of singular vectors. The orthogonality property usually degrades as one estimates more singular values. Thus, \cite{roy2021new} propose to use a Gram-Schmidt orthogonalization trick~\citep{giraud2005loss} in between the iterative equations. In particular, between alternatively using~\eqref{eqn:general-1}-\eqref{eqn:general-3}, the estimates of the $k$-th singular vectors $\bb{a}_k$ and $\bb{b}_k$ are updated as $\bb{a}_k \leftarrow \bb{a}_k - \sum_{r = 1}^{(k-1)} {\bb{a}_k}\tr \bb{a}_r$ and $\bb{b}_k \leftarrow \bb{b}_k - \sum_{r = 1}^{(k-1)} {\bb{b}_k}\tr \bb{b}_r$ for all but the first singular value. Further details about the estimation algorithm can be found in~\cite{roy2021new}.

\section{Mathematical Properties}\label{sec:math-property}

Since the estimation process of the subsequent singular values and vectors follows from the same estimation technique of rank one decomposition on the residual matrix, we shall restrict our attention to the study of the properties of the rSVDdpd estimator for the rank one decomposition only. In order to make sure all the results developed in this section are true for the subsequent singular values and vectors as well, the assumptions on the distribution of the data matrix $\bb{X}$ must hold for the residual matrix after subtracting the effect of previous singular values and vectors. We assume that this is the case. This is indeed true for the situations where the true distributions of the random variable $X_{ij}$ denoting the elements of the data matrix $\bb{X}$ belong to a location family of distributions with location parameters $\sum_{r = 1}^k \lambda_r a_{ir} b_{jr}$, where the vectors $\bb{a}_1, \ldots \bb{a}_r$ and $\bb{b}_1, \ldots \bb{b}_r$ are of unit $L_2$ norm and $\lambda_r$ denotes the true singular values of the matrix.

Before proceeding with the description of the mathematical properties of the rSVDdpd estimator, we fix some notations to be used throughout the paper. For a matrix $\bb{X}$, its entries will be denoted by $X_{ij}$. For any parameter $\bb{\theta}$, the superscript $\bb{\theta}^g$ will denote the true value of the parameter that is being estimated, when the true data density is $g$, the superscript $\bb{\theta}^{(t)}$ will denote the value of the estimated parameter at $t$-th iteration of the algorithm and $\bb{\theta}^\ast$ will be used to indicate its limit, provided that the sequence of iterated estimates converge. In the asymptotic analysis, we use the notation $a_n = \mathcal{O}(b_n)$ to denote the scenario when for sufficiently large $n$, there exists constant $C$, independent of $n$ such that $\vert a_n\vert \leq C\vert b_n\vert$ for any two sequences $\{ a_n\}$ and $\{ b_n\}$. On the other hand, the notation $a_n \asymp b_n$ is used for asymptotic equivalence, i.e., there exists two constants $0 < C_1 < C_2 < \infty$ such that $C_1a_n < b_n < C_2a_n$ for sufficiently large $n$.

In the context of rank-one estimation, we redefine the parameter as $\bb{\theta} = \left( \lambda, \{ u_i \}_{i=1}^n, \{ v_j\}_{j=1}^p, \sigma^2 \right)$, comprising the first singular value and corresponding vectors. The corresponding parameter space becomes $\bb{\Theta} = (0,\infty) \times S_n^+ \times S_p \times [\epsilon, \infty)$, where $\epsilon > 0$ is a small positive quantity and $S_n^+$ is the part of $n$-dimensional unit hyperspheres centered at the origin and lying in the positive orthant, i.e.
\begin{equation}
    S_n^+ = \left\{ (x_1, \ldots x_r) : \sum_{i=1}^r x_i^2 = 1, x_i \geq 0, i = 1, \ldots r \right\}\label{eqn:sr-plus-defn}
\end{equation}
\noindent Such a restriction on the parameter space is required to ensure identifiability of the robust SVD problem, since one can switch the signs of both $\bb{u}$ and $\bb{v}$ resulting in the same decomposition. In view of~\eqref{eqn:general-1}-\eqref{eqn:general-3}, the sequence of estimates $\bb{\theta}^{(t)}$ is then related in the following way,
\begin{align}
    \lambda^{(t+1)} u_i^{{(t+1)}} & = \dfrac{\sum_j v_j^{(t)} X_{ij} \psi(e_{ij}^{(t)}/\sigma^{(t)}) }{\sum_j (v_j^{(t)})^2 \psi(e_{ij}^{(t)}/\sigma^{(t)}) }
    \label{eqn:iteration-rule-1}  \\
    \lambda^{(t+1)} v_j^{{(t+1)}} & = \dfrac{\sum_i u_i^{(t+1)} X_{ij} \psi(e_{ij}^{(t)}/\sigma^{(t)}) }{\sum_i (u_i^{(t+1)})^2 \psi(e_{ij}^{(t)}/\sigma^{(t)}) }
    \label{eqn:iteration-rule-2}  \\
    (\sigma^2)^{(t+1)} = & \dfrac{(np)^{-1}\sum_{i,j} (e_{ij}^{(t)})^2 \psi(e_{ij}^{(t)}/\sigma^{(t)}) }{ (np)^{-1}\sum_{i,j} \psi(e_{ij}^{(t)}/\sigma^{(t)}) - \frac{\alpha}{1+\alpha} M_f}
    \label{eqn:iteration-rule-3}
\end{align}
\noindent for all $t = 0, 1, \dots$. Before proceeding with the statistical properties of the rSVDdpd estimator, we establish the convergence of the above iterative procedure under two simple assumptions.
\begin{enumerate}[label = (A\arabic*), ref = (A\arabic*)]
    \item\label{assum:f-diff} Assume that the model density $f$ is twice differentiable with respect to its arguments.
    \item\label{assum:f-diff-cond} The model density $f$ is symmetric and satisfies $f'(x) \leq 0$ for all $x > 0$ and
          \begin{equation}
              \dfrac{1}{x} > \alpha \dfrac{f'(x)}{f(x)} + u'(x) \dfrac{f(x)}{f'(x)}, \ \text{for } x > 0.
              \label{eqn:psi-prime-cond}
          \end{equation}
    \item\label{assum:psi-bound} There exists a constant $K$ such that $x^2 \psi(x) < K$ for all $x \geq 0$.
\end{enumerate}

\begin{remark}
    The Assumption~\ref{assum:f-diff-cond} does not impose strict conditions on the choice of $f$. As a result of $f'(x) \leq 0$, we obtain $\psi(x) \geq 0$, and the condition~\eqref{eqn:psi-prime-cond} implies that $\psi$ is decreasing. Thus, the provided conditions simply amount to the requirement that the weights in~\eqref{eqn:general-1}-\eqref{eqn:general-3} are nonnegative, and the larger errors are down-weighted so that it leads to a robust estimator.
\end{remark}

\begin{remark}
    When $f$ is standard normal density, $\psi(x) = e^{-\alpha x^2/2}$. This ensures that Assumption~\ref{assum:f-diff-cond} is satisfied due to the nonnegativity and decreasing nature of $\psi$. Assumption~\ref{assum:psi-bound} is satisfied in this case with $K = 2/\alpha e$.
\end{remark}

\begin{theorem}\label{thm:convergence}
    For a fixed $n$ and $p$ and the data matrix $\bb{X}$, if assumptions~\ref{assum:f-diff} and~\ref{assum:f-diff-cond} hold, then the sequence of estimates $\bb{\theta}^{(t)}$ obtained through~\eqref{eqn:iteration-rule-1}-\eqref{eqn:iteration-rule-3} converges to a local minimizer of $H_\alpha^{(1)}(\bb{\theta})$ shown in Eq.~\eqref{eqn:H-general}.
\end{theorem}

Let us denote this converged rSVDdpd estimator as $\bb{\theta}^\ast = (\lambda^\ast, \{ u_i^\ast \}_{i = 1}^{n}, \{ v_j^\ast \}_{j=1}^{p}, (\sigma^\ast)^2 )$. On the other hand, let us denote the population counterpart as $\bb{\theta}^g = (\lambda^g, \{ u_i^g \}_{i = 1}^{n}, \{ v_j^g \}_{j=1}^{p}, (\sigma^g)^2 )$, which is the true value of the parameters that are ultimately being estimated. Similar to the iteration rules~\eqref{eqn:iteration-rule-1}-\eqref{eqn:iteration-rule-3} for the rSVDdpd estimator, the true value $\bb{\theta}^g$ is also expected to satisfy such fixed point criteria, in the sense of overall population based measures rather than its empirical counterparts. With this in mind, we start by defining these ``best'' fitting parameters for the particular setup of SVD under consideration.

\begin{definition}
    Let, $\bb{X}$ be a data matrix of order $n \times p$ such that its $(i,j)$-th entry $X_{ij}$ follows a distribution with density function $g_{ij}$, for all $i = 1, \ldots n$ and $j = 1, \ldots p$ and $X_{ij}$s are independent to each other. Then, $\bb{\theta}^g = (\lambda^g, \{ u_i^g \}_{i = 1}^{n}, \{ v_j^g \}_{j=1}^{p}, (\sigma^g)^2 )$ is called a ``best'' fitting parameter if the following conditions hold
    \begin{enumerate}
        \item The $\{u_i^g \}$s and $\{ v_j^g \}$s constitute entries of unit vectors, i.e.,
              \begin{equation}
                  \sum_{i=1}^{n} (u_i^g)^2 = 1, \ \text{ and } \
                  \sum_{j=1}^{p} (v_j^g)^2 = 1.
                  \label{eqn:true-norm}
              \end{equation}
        \item For any $i = 1, 2, \dots n,\ j = 1, 2, \dots p$,
              \begin{equation}
                  \lambda^g u_i^g = \argmin_{a} \int V_f(\cdot; a, v_j^g, (\sigma^g)^2 ) g_{ij}.
                  \label{eqn:true-ai}
              \end{equation}
        \item For any $i = 1, 2, \dots n,\ j = 1, 2, \dots p$,
              \begin{equation}
                  \lambda^g v_j^g = \argmin_{b} \int V_f(\cdot; u_i^g, b, (\sigma^g)^2 ) g_{ij}.
                  \label{eqn:true-bj}
              \end{equation}
        \item For any $i = 1, 2, \dots n,\ j = 1, 2, \dots p$,
              \begin{equation}
                  (\sigma^g)^2 = \argmin_{\sigma^2}  \int V_f(\cdot; \lambda^g u_i^g, v_j^g, \sigma^2) g_{ij} \\
                  = \argmin_{\sigma^2}  \int V_f(\cdot; u_i^g, \lambda^g v_j^g, \sigma^2) g_{ij}.
                  \label{eqn:true-sigma}
              \end{equation}
    \end{enumerate}
    \label{defn:best-parameter}
    \noindent where
    \begin{equation}
        V_f(x; a, b, \sigma^2) \\
        = \sigma^{-\alpha} \left[ M_f - \left(1 + \dfrac{1}{\alpha} \right) f^\alpha\left( \left\vert \dfrac{x - ab}{\sigma} \right\vert \right) \right].
        \label{eqn:v-ab-sigma}
    \end{equation}
\end{definition}
\noindent Here,~\eqref{eqn:true-ai} shows that the minimizer of the quantity on the right-hand side of the equation is always $\lambda^g a_i^g$, independent of the choice of column index $j$. This Assumption holds if the true densities $g_{ij}$ are densities of the normal distributions with the location parameters being elements from the best rank one approximation of $\bb{X}$, i.e. the entries of the data matrix $X_{ij}$s are normally distributed with mean $\mu_{ij}$ and constant variance $\sigma^2$, and the matrix $\bb{\mu} = \left( \mu_{ij} \right)_{i=1, j=1}^{n,p}$ is of unit rank.

\subsection{Uniqueness}

Since the aim of the rSVDdpd algorithm is to robustly estimate the singular values and the singular vectors of a given data matrix, it is required to show that the ``best'' fitting parameters as introduced by Definition~\ref{defn:best-parameter} resemble the behaviour of the usual singular values and vectors. Regarding this, the following theorem claims that if the elements of the data matrix $\bb{X}$ follows a decomposition as in~\eqref{eqn:rsvd-decomp}, then the ``best'' fitting parameter given by Definition~\ref{defn:best-parameter} matches with the usual singular values and vectors.
\begin{theorem}
    \label{thm:def-match}
    Let the data matrix $\bb{X}$ be such that $X_{ij} = \lambda^\ast u_i^\ast v_j^\ast + \epsilon_{ij}$ where $\epsilon_{ij}$s are i.i.d. random variables with density $(\sigma^\ast)^{-1}f(\cdot/\sigma^\ast)$. Then $\bb{\theta}^g$ is the unique ``best'' fitting parameter if $\bb{\theta}^g = \left( \lambda^\ast, \{u_i^\ast \}_{i=1}^n, \{ v_j^\ast \}_{j=1}^p, (\sigma^\ast)^2  \right)$ belongs to the parameter space $\bb{\Theta}$.
\end{theorem}

The following corollaries of Theorem~\ref{thm:def-match} are now immediate which establish the validity of the best fitting parameter, for the case when the entries of the data matrix $\bb{X}$ follow a normal distribution or are deterministic (which is a special case of the family of normal distribution with variance parameter equal to $0$).
\begin{corollary}
    \label{cor:normal-match}
    Let the data matrix $\bb{X}$ be such that $X_{ij}$ are i.i.d. $N(\lambda^\ast u_i^\ast v_j^\ast, (\sigma^\ast)^2)$. Then $\bb{\theta}^g$ is the unique ``best'' fitting parameter with normal model family of densities if $\bb{\theta}^g = \left( \lambda^\ast, \{u_i^\ast \}_{i=1}^n, \{ v_j^\ast \}_{j=1}^p, (\sigma^\ast)^2  \right) \in \bb{\Theta}$. In this case, the $V_f$-function as in~\eqref{eqn:v-ab-sigma} is denoted by $V_\phi$ and is given as
    \begin{equation}
        V_\phi(x; a, b, \sigma^2)
        = \sigma^{-\alpha} \left( \dfrac{1}{\sqrt{1+\alpha}} - \left(1+\dfrac{1}{\alpha} \right) \phi^\alpha((x-ab)/\sigma) \right)
        \label{eqn:V-normal-function}
    \end{equation}
    \noindent where $\phi$ is the standard normal density function.
\end{corollary}

An immediate corollary of Theorem~\ref{cor:normal-match} for a deterministic data matrix $\bb{X}$ follows from the observation that degenerate distribution is one special case of the family of normal distributions with the variance parameter being equal to $0$.

\begin{corollary}
    \label{cor:degenerate-match}
    If the data matrix $\bb{X}$ is of rank $1$ such that $X_{ij} = \lambda^\ast u_i^\ast v_j^\ast$ for all $i = 1, 2, \dots n,\ j = 1, 2, \dots p$, with $\sum_i (u_i^\ast)^2 = \sum_j (v_j^\ast)^2 = 1$, then there exists a unique ``best'' fitting parameter given by $\bb{\theta}^g = (\lambda^\ast, \{ u_i^\ast \}_{i = 1}^{n}, \{ v_j^\ast \}_{j=1}^{p}, 0 )$ if $\bb{\theta}^g \in \bb{\Theta}$.
\end{corollary}
\noindent In the deterministic setup given in Corollary~\ref{cor:degenerate-match}, since the true distribution is a degenerate distribution at $x_{ij} = \lambda^\ast a_i^\ast b_j^\ast$, it follows that $\int V(x; c, d, \sigma^2) g_{ij}(x)dx = V(\lambda^\ast a_i^\ast b_j^\ast; c, d, \sigma^2)$. Therefore, the population version of the alternating estimating equations given in~\cite{roy2021new} becomes the same as the iteration formulas given in~\eqref{eqn:iteration-rule-1}-\eqref{eqn:iteration-rule-3}. In other words, $\bb{\theta}^g = (\lambda^\ast, \{ u_i^\ast \}_{i = 1}^{n}, \{ v_j^\ast \}_{j=1}^{p}, 0 )$ becomes the unique fixed point of the alternating iteration rules, in the restricted parameter space $\bb{\Theta}$.

\subsection{Equivariance Properties}

We have seen that if the entries of $\bb{X}$ are of special structure, i.e., low-rank plus i.i.d. noise, under correct specification of the family of densities $f$, the best fitting parameter is equivalent to the singular values and corresponding vectors. However, in case of misspecification or if $\bb{X}$ do not follow the specific structure, but has independent entries, we can show that a ``best'' fitting parameter satisfies equivariance properties similar to the singular values and vectors. One such equivariance property is that whenever the data matrix is multiplied by some scalar quantity, the singular values are also multiplied by the same scalar. However, the singular vectors remain unchanged. The following theorem presents the equivariance property for a ``best'' fitting parameter.
\begin{theorem}
    \label{thm:scale-equivariance}
    If a ``best'' fitting parameter for matrix $\bb{X}$ is $\bb{\theta}^g$, then a ``best'' fitting parameter for the matrix $c\bb{X}$ is $\widetilde{\bb{\theta}}^g = (c\lambda^g, \{ a_i^g \}_{i = 1}^{n}, \{ b_j^g \}_{j=1}^{p}, (c\sigma^g)^2 )$ for any real constant $c$.
\end{theorem}

Another equivariance property of the usual singular value decomposition is that under any row (or column) permutation of the data matrix, the entries of the left (or right) singular vectors permute accordingly, while the singular values remain unaffected. Such a property also holds for a ``best'' fitting parameter. Let $\pi_R$ and $\pi_C$ denote some such permutation on the row and column indices of the matrix respectively.
\begin{theorem}
    \label{thm:permute-equivariance}
    Let $\bb{P}, \bb{Q}$ are permutation matrices corresponding to the permutations $\pi_R$ and $\pi_C$, on the row and column indices of the matrix $\bb{X}$. If a ``best'' fitting parameter for the matrix $\bb{X}$ is $\bb{\theta}^g$, then a ``best'' fitting parameter for permuted matrix $\bb{P}\bb{X}\bb{Q}\tr$ is the correspondingly permuted version of $\bb{\theta}^g$ given by $\widetilde{\bb{\theta}}^g = (\lambda^g, \{ u_{\pi_R(i)}^g \}_{i = 1}^{n}, \{ v_{\pi_C(j)}^g \}_{j=1}^{p}, (\sigma^g)^2  )$.
\end{theorem}

While Theorem~\ref{thm:scale-equivariance} and Theorem~\ref{thm:permute-equivariance} indicate the connection between singular values and a ``best'' fitting parameter, similar equivariance property is also obeyed by the converged rSVDdpd estimator. The following theorems demonstrate the same.
\begin{theorem}
    \label{thm:scale-equivariance-estimate}
    Let $\bb{\theta}^\ast$ be the converged rSVDdpd estimator for the matrix $\bb{X}$, starting from an initial estimate $\bb{\theta}^{(0)}$. Then, for any constant $c \in \R \backslash \{0\}$, the rSVDdpd estimator for the data matrix $\bb{X}\dash = c\bb{X}$ converges to $(c\lambda^\ast, \{ u_i^\ast \}_{i = 1}^{n}, \{ v_j^\ast \}_{j=1}^{p}, (c\sigma^\ast)^2 )$ provided the initial estimate is $(c\lambda^{(0)}, \{ u_i^{(0)} \}_{i = 1}^{n}, \{ v_j^{(0)} \}_{j=1}^{p}, (c\sigma^{(0)})^2 )$.
\end{theorem}

\begin{theorem}
    \label{thm:permute-equivariance-estimate}
    Let $\bb{\theta}^\ast$ be the converged rSVDdpd estimator for the matrix $\bb{X}$, starting from an initial estimate $\bb{\theta}^{(0)}$. Also, let $\bb{P}$ and $\bb{Q}$ be the permutation matrices corresponding to the permutations $\pi_R$ and $\pi_C$ respectively. Then, starting with the new initial estimate $(\lambda^{(0)}, \{ u_{\pi_R(i)}^{(0)} \}_{i = 1}^{n}, \{ v_{\pi_C(j)}^{(0)} \}_{j=1}^{p}, (\sigma^{(0)})^2 )$, the rSVDdpd estimator for the data matrix $\bb{PXQ}\tr$ converges to the corresponding permuted version of $\bb{\theta}^\ast$ given by $(\lambda^\ast, \{ u_{\pi_R(i)}^\ast \}_{i = 1}^{n}, \{ v_{\pi_C(j)}^\ast \}_{j=1}^{p}, (\sigma^\ast)^2 )$.
\end{theorem}

\subsection{Consistency of the rSVDdpd estimator}\label{sec:consistency-thm}

The convergence theorem presented in Theorem~\ref{thm:convergence} ensures that under some minimal assumptions, the rSVDdpd algorithm given by the iterations~\eqref{eqn:iteration-rule-1}-\eqref{eqn:iteration-rule-3} converges to the rSVDdpd estimator, i.e., a local minimizer $\bbhat{\theta}^\ast$ of $H_\alpha^{(1)}(\bb{\theta})$ given in~\eqref{eqn:H-general}. However, in view of Definition~\ref{defn:best-parameter} of a ``best'' fitting parameter $\bb{\theta}^g$, it is necessary to know whether such a local optimum remains close to a ``best'' fitting parameter in an asymptotic sense. In this asymptotic regime, we allow both the matrix dimensions $n$ and $p$ to grow to infinity, subject to a constant ratio in limit, i.e., $n/p \rightarrow c$ for some $c \in (0, \infty)$.

Answering this question about statistical consistency of the rSVDdpd estimator has two technical difficulties. Firsly, the parameter space $\bb{\Theta}$ is not necessarily convex due to the presence of the coordinates related to singular vectors. The first problem can be circumvented using an inverse stereographic projection which transforms this non-convex parameter space $\bb{\Theta}$ into a convex parameter space $\bb{\Xi} \subseteq \R^{(n+p)}$. We call this parameter space $\bb{\Xi}$ as the natural parameter space in the given setup. The one-one transformation $\mathcal{T}$ between these two parameter spaces $\bb{\Theta}$ and $\bb{\Xi}$ are governed by the following two equations
\begin{equation}
    \mathcal{T}(\lambda, \{u_i\}_{i=1}^{n}, \{ v_j\}_{j=1}^p, \sigma^2)
    = \left( \lambda, \left\{ \frac{u_i}{(1-u_n)} \right\}_{i=1}^{(n-1)}, \left\{ \frac{v_j}{(1-v_n)} \right\}_{j=1}^{(p-1)}, \sigma^2 \right),
    \label{eqn:stereo-theta-to-eta}
\end{equation}
\noindent and,
\begin{equation}
    \mathcal{T}^{(-1)}\left( \lambda, \{\alpha_i\}_{i=1}^{(n-1)}, \{ \beta_j\}_{j=1}^{(p-1)}, \sigma^2  \right) \\
    = \left( \lambda, \left\{ \frac{2\alpha_i}{U^2 + 1} \right\}_{i=1}^{(n-1)}, \frac{U^2 - 1}{U^2+ 1},  \left\{ \frac{2\beta_j}{V^2 + 1} \right\}_{j=1}^{(p-1)}, \frac{V^2 - 1}{V^2+ 1}, \sigma^2 \right),
    \label{eqn:stereo-eta-to-theta}
\end{equation}
\noindent where $U^2 = \sum_{i=1}^{(n-1)} \alpha_i^2$ and $V^2 = \sum_{j=1}^{(p-1)} \beta_j^2$. Accordingly, we denote $\bb{\eta}$ as an element of this natural parameter space $\bb{\Xi}$, where the corresponding transformed parameter $\bb{\theta} = \mathcal{T}^{(-1)}(\bb{\eta})$ denotes an element of $\bb{\Theta}$.

The second problem is that the length of the singular vectors is not fixed and grows with the dimension of the matrix. Thus, as $n, p\rightarrow \infty$, the dimension of the parameter space, i.e., $(n + p + 2)$ grows linearly in $n$ or $p$ and increases to infinity. This means, although, rSVDdpd uses the MDPDE proposed by~\cite{ghosh2013} to solve the linear regression problems, their consistency results cannot be used directly as it assumes the dimension of the parameter space to be fixed. This varying dimension problem has been of considerable interest to many authors under the M-estimation setup~\citep{huber1973robust,portnoy1984asymptotic}. Most of these results assume convexity of the objective function~\citep{he2000parameters}, but that cannot be employed in our case. Here, the objective function is convex in any of the parameters individually when the other parameters are kept fixed, but becomes non-convex if all the parameters are taken together. To deal with this problem, we need to apply some concentration bounds on the error terms, for which we will restrict our attention to the scenario where the model density $f$ is the standard normal density function.

Now that the necessary foundations are laid out, we can present the consistency theorem, which claims that under some reasonable assumptions indicatd below, the minimizer $\bb{\theta}^\ast$ of $\mathcal{H}^{(1)}_\alpha$ as given in Theorem~\ref{thm:convergence}, is a consistent estimator of the best fitting parameter $\bb{\theta}^g$. However, note that the description of a ``best'' fitting parameter indicated in Definition~\ref{defn:best-parameter} is applicable for fixed $n$ and $p$. In contrast, statistical consistency is an asymptotic property, which requires the dimensions of the matrix $n$ and $p$ to tend towards infinity. To resolve this conflict in a unified setup, we assume that a sequence of ``best'' fitting parameters exists for each fixed $n$ and $p$. We denote the $(i, j)$-th entry of the data matrix $\bb{X}_{n,p}$ of order $n \times p$ by the random variable $(X_{ij})_{n,p}$, and make several assumptions about it as follows
\begin{enumerate}[label = (B\arabic*), ref = (B\arabic*)]
    \item\label{assum:best-fitting-distn} There exists a sequence $\bb{\theta}^g_{n, p} = (\lambda^g, \{ u_{i, n}^g \}_{i=1}^n, \{ v_{j, p}^g \}_{j=1}^p, (\sigma^g_{n,p})^2)$ of the best fitting parameters such that $X_{ij}$s are independently distributed as
  \begin{equation*}
      X_{ij} = \lambda^g u_{i,n}^g v_{j,p}^g + \sigma^g_{n,p} Z_{ij},
  \end{equation*}
  \noindent for all $i = 1, 2, \dots n; j = 1, 2, \dots p$ and $Z_{ij}$s are i.i.d. following a common density $g$.
    \item \label{assum:dens-g-4th} The density function $g$ is such that the integrals $\int e^{-\alpha z^2/2} g(z)dz$ exist, and are three times differentiable and the derivatives can be taken under the integral sign. Also, the integrals $\int z^k e^{-\alpha z^2/2} g(z)$ are finite for $k = 0, 1, \dots 4$.
    \item\label{assum:dens-g-symmetric} The density function $g$ is symmetric.
    \item\label{assum:open-rect} For each pair of positive integers $n$ and $p$, there exists an open rectangle $S_{n,p}$ inside $\bb{\Xi}_{n,p} \subset \R^{(n+p)}$ containing the natural parametrization $\bb{\eta}^g_{n,p}$ of $\bb{\theta}^g_{n,p}$ such that the sequence of sets $\mathcal{\tau}(S_{n,p})$ does not have a limit point with $\bb{a} = (0, \dots, 0, 1)$ or $\bb{b} = (0, \dots, 0, 1)$.
    \item\label{assum:stereographic} The converged rSVDdpd estimate for the data matrix $\bb{X}_{n,p}$, i.e, $\bb{\theta}_{n,p}^\ast$ (the minimizer of $\mathcal{H}_{n,p}$ as indicated by Theorem~\ref{thm:convergence}) satisfy $(\bb{a}_{n}^{\ast})\tr \neq (0, \dots, 0, 1)$ and $(\bb{b}_p^{\ast})\tr \neq (0, \dots, 0, 1)$ for all sufficiently large $n$ and $p$.
    \item\label{assum:variance-reduce} The variance $(\sigma^g_{n,p})^2$ in best fitting parameters satisfy $(\sigma^g_{n,p})^2 \asymp (np)^{-1/2}$.
\end{enumerate}
The first Assumption~\ref{assum:best-fitting-distn} is simply a description of the setup. Assumptions~\ref{assum:dens-g-4th} and~\ref{assum:open-rect} are standard assumptions connected to the MDPDE~\citep{ghosh2013}. Assumptions~\ref{assum:dens-g-symmetric} and~\ref{assum:variance-reduce} are required to provide concentration bound on the covariance terms and tail probabilities respectively. It is well known from random matrix theory that the Gaussian ensemble with each entry following a standard normal distribution has the singular values asymptotically at the order of $(\sqrt{n} + \sqrt{p})$~\citep{tracy1993introduction,mehta2004random}. However, since Assumption~\ref{assum:best-fitting-distn} indicates that the same $\lambda^g$ acts as the singular value for any $n$ and $p$, the variance of the entries of the data matrix has to go down asymptotically to ensure that the singular value does not grow with increasing $n$ or $p$.

\begin{theorem}
    Under Assumptions~\ref{assum:best-fitting-distn}-\ref{assum:variance-reduce}, if the model density $f$ is the standard normal density, then there exists a sequence of rSVDdpd estimates $\bb{\theta}_{n,p}^\ast$ which is consistent to a sequence of ``best'' fitting parameters $\bb{\theta}_{n,p}^g$. That means, as both dimensions of the data matrix $\bb{X}_{n,p}$ (i.e. $n$ and $p$) tend to infinity subject to a constant ratio in limit, i.e. $\lim_{\substack{n\rightarrow \infty\\ p \rightarrow \infty}} \frac{n}{p} = c$ for some $c \in (0, \infty)$,
    $\norm{\bb{\theta}_{n,p}^\ast - \bb{\theta}_{n,p}^g} \rightarrow 0$ in probability.
    \label{thm:consistency}
\end{theorem}

Since \cite{roy2021new} derived rSVDdpd as an extension of MDPDE in linear regression setup based on the work of~\cite{ghosh2013}, it is natural to think that the result on the consistency of the MDPDE given in Theorem 3.1 of the same paper can be imitated to deliver a proof of the consistency of the proposed rSVDdpd estimators. However, there are several major complications involved.
\begin{enumerate}
    \item The basic assumption required for consistency of the MDPDE in the INH (independent and non-homogeneous) setup considered in~\cite{ghosh2013} is the existence of an open set in the parameter space. However, in this particular setup, the parameter space $\bb{\Theta}$ by itself does not contain any open neighbourhood. Therefore, all the necessary formulations are required to be applied on the natural parametrization $\bb{\eta} \in \bb{\Xi}$ instead, which converts the setup into a nonlinear regression problem.
    \item In the case of SVD, the length of the singular vectors is not fixed and grows with the dimension of the matrix. Thus, as $n, p\rightarrow \infty$, the dimension of the parameter space also increases to infinity. This problem has been of considerable interest to many authors under the M-estimation setup~\citep{huber1973robust,portnoy1984asymptotic}. All of these results assume convexity of the objective function~\citep{he2000parameters}, but that cannot be employed in our case. Here, the objective function is convex in any of the parameters individually when the other parameters are kept fixed, but becomes non-convex if all the parameters are taken together.        
    \item In each of the alternating iterations, the consistency ensures that the resulting estimator based on the minimization of that particular iteration is probabilistically close to the minimizer of the population version of the iterative equation. However, the population version of the iterative equation depends on the current estimates of the other parameters. Hence, in each of the iterations, the empirical estimates are allowed to deviate from the population estimates and these errors sum up as the number of iterations increases. Hence, we require a bound on the tail of the distribution of such a sum of errors to ensure the consistency of the rSVDdpd estimator.
\end{enumerate}

\noindent Thus, a non-trivial modification of the existing proof technique is required. Due to its length and complications, the proof of this theorem is deferred to Appendix~\ref{appendix:consistency-proof}.

\begin{remark}
    Theorem~\ref{thm:scale-equivariance-estimate} and~\ref{thm:permute-equivariance-estimate} are, respectively, the empirical counterparts of Theorem~\ref{thm:scale-equivariance} and~\ref{thm:permute-equivariance}. In view of Theorem~\ref{thm:scale-equivariance-estimate} and~\ref{thm:permute-equivariance-estimate}, the equivariance properties hold given that the initial values of the iterations of rSVDdpd also satisfy such equivariance properties. However, from the convergence and the consistency of the rSVDdpd estimator, it follows that for large $n$ and $p$, the converged estimator can be made arbitrarily close to the true ``best'' fitting parameters. Since these ``best'' fitting parameters obey equivariance properties as assured in Theorem~\ref{thm:scale-equivariance} and~\ref{thm:permute-equivariance}, it follows that for large $n$ and $p$, the converged estimator will also approximately satisfy these equivariance properties, irrespective of the equivariance of starting values.
\end{remark}

Instead of restricting the rSVDdpd estimator to only the normal family of model densities, one can take $f$ to be any subgaussian density. In this case, Assumption~\ref{assum:dens-g-4th} needs to be appropriately modified to ensure that the corresponding $\psi(\cdot)$ function is two-times continuously differentiable and $\int z^k \psi(z) g(z)dz$, $\int z^2(\psi'(z))^2 g(z)dz$ and $\int  z^2 \psi''(z) g(z)$, are all bounded.

Also note that Theorem~\ref{thm:consistency} ensures the consistency of the rSVDdpd under a general setup where the errors follow any arbitrary density function $g$ subject to the Assumptions~\ref{assum:best-fitting-distn}-\ref{assum:variance-reduce}. Therefore, it also allows density functions of the form $g = (1 - \epsilon)g_1 + \epsilon g_2$, which is $\epsilon$-contaminated version of density $g_1$ contaminated by density $g_2$, provided that both $g_1$ and $g_2$ are symmetric functions. Additionally, to ensure that Assumption~\ref{assum:dens-g-4th} is satisfied, a sufficient condition is that the density function $g$ is thrice continuously differentiable and the random variable $Z$ with density $g$ has finite fourth-order moments. However, even if the moment condition does not hold, one can directly show Assumption~\ref{assum:dens-g-4th} for all $\alpha > 0$. For instance, in the case of Cauchy density,
\begin{equation*}
    \E(Z^4 e^{-\alpha Z^2/2}) = \int_{-\infty}^{\infty} \dfrac{x^4e^{-\alpha x^2/2}}{\pi (1 +x^2)} = \dfrac{2}{\sqrt{\pi}}e^{\alpha/2} \int_{\sqrt{\alpha}/\sqrt{2}}^\infty e^{-t^2}dt - \dfrac{\alpha - 1}{\sqrt{2\pi} \alpha^{3/2}},
\end{equation*}
\noindent for all $\alpha > 0$, and is finite. Therefore, except for $\alpha = 0$ (i.e., MLE), the consistency of the rSVDdpd estimator is ensured when errors follow the Cauchy distribution with heavy tails.

\section{Numerical Illustrations: Simulation Studies}\label{sec:simulation}

In this section, we compare the performance of the rSVDdpd estimator with two existing robust SVD estimators, namely the ones proposed by~\cite{niss2001} and~\cite{zhang2013}. Implementation of the robust SVD algorithm proposed by~\cite{niss2001} (to be referred to here as pcaSVD) is available as an R package pcaMethods~\citep{pcaMethods}, which outputs all singular values and vectors of the input data matrix. The second algorithm by \cite{zhang2013} obtains the first pair of singular vectors based on the minimization procedure
\begin{equation*}
(\widehat{\bb{u}}, \widehat{\bb{v}}) = \argmin_{(\bb{u}, \bb{v})} \left[ \rho\left( \dfrac{\bb{X} - \bb{u}\bb{v}\tr}{\sigma} \right) + \mathcal{P}_\lambda(\bb{u}, \bb{v})\right]    
\end{equation*}
\noindent where $\rho(\cdot)$ is a robust loss function (namely Huber's loss function) and $\mathcal{P}_\lambda$ is a regularization penalty term to motivate smoothing in the entries of the singular vectors. For an extensive comparison, we consider two variants of this algorithm. In one variation, we perform the minimization with only Huber's loss function without any penalty term, which we shall refer to as RobSVD, while in the other variation, we follow the recommended procedure of minimization with penalty term, which we shall call RobRSVD. Implementation of both of these variants is available in the R package RobRSVD~\citep{RobRSVD} which is programmed only to output the first singular value and its corresponding singular vectors. Thus, in order to compare the performances of the algorithms on an equal footing, we add a wrapper outside this function to apply the same algorithm on the residual matrix to output the subsequent singular values (see~\cite{niss2001} for details). Along with these, we also consider an implementation of the usual $L_2$ norm minimization based SVD procedure available in R base package written using LAPACK Fortran library. For the rSVDdpd algorithm, we use the R package ``rsvddpd'' as provided by~\cite{roy2021new}.

\subsection{Simulation Setups}\label{sec:sim-setup}

To compare the performance of the robust SVD approaches, we employ a Monte Carlo method. We generate random errors and add them to the true matrix, then apply each of the robust SVD algorithms and track the estimates obtained from each sample. The mean squared error (MSE) and bias are computed for each singular value based on $B = 1000$ Monte Carlo samples, serving as the accuracy measures. The sum of squared biases and the sum of MSE across all singular values are calculated for each estimator.

For comparing the left and right singular vectors, we consider a dissimilarity score between two normalized vectors, denoted as $\text{Diss}(\bb{u}, \bb{v})$, which is equal to $1 - \vert \langle \bb{u}, \bb{v}\rangle \vert$. Here, $\vert x\vert$ represents the absolute value of $x$, and $\langle \bb{u}, \bb{v}\rangle$ is the Euclidean inner product between two vectors $\bb{u}$ and $\bb{v}$. The dissimilarity score is $0$ if $\bb{u} = \bb{v}$ or $\bb{u} = (-\bb{v})$, and it is equal to $1$ if $\bb{u}$ and $\bb{v}$ are orthogonal. The average dissimilarity score between the estimated singular vectors and the true singular vectors, computed over all Monte Carlo samples, is used as a performance measure.
 
To construct the true data matrix $\bb{X}$ for a given singular value decomposition, we use the coefficients of the first three orthogonal polynomial contrasts of order $10$ and $4$ and arrange them as columns of matrices $\bb{U}$ and $\bb{V}$ respectively. The resulting data matrix $\bb{X}$ is then constructed as 
\begin{equation}
    \bb{X} = \bb{U} \begin{bmatrix}
        10 & 0 & 0\\
        0 & 5 & 0\\
        0 & 0 & 3\\ 
     \end{bmatrix} \bb{V}\tr,
     \label{eqn:sim-X-matrix}  
\end{equation}
\noindent with the true singular values being $10, 5$ and $3$. In each of the resamples, errors following a pre-specified distribution are added to the entries of $\bb{X}$. Based on the chosen distribution, we divide the simulation scenarios broadly into $5$ categories denoted by~\ref{sim:setup-1}-\ref{sim:setup-5}, as follows.
\begin{enumerate}[label = (S\arabic*), ref = (S\arabic*)]
    \item\label{sim:setup-1} The errors follow the standard normal distribution $\normdist(0, 1)$ (no outliers per se).
    \item\label{sim:setup-2} The errors are distributed according to a contaminated standard normal distribution namely
    $$e_{ij} \sim (1 -\epsilon) \normdist(0,1) + \epsilon \delta_{25},$$
    where $\epsilon$ is the amount of contamination and $\delta_{25}$ denotes the degenerate distribution at $25$. Based on the amount of contamination, we consider three sub-cases of this simulation setting.
    \begin{enumerate}[label = (S2\alph*), ref = (S2\alph*)]
        \item\label{sim:setup-2a} $\epsilon = 0.05$, denotes only $5\%$ contamination, which corresponds to a relatively light amount of outlying observations.
        \item\label{sim:setup-2b} $\epsilon = 0.1$, denoting medium level contamination with the presence of $10\%$ outlying values.
        \item\label{sim:setup-2c} $\epsilon = 0.2$, denoting heavy contamination with approximately $20\%$ of the entries being same as the outlying observation of $25$.
    \end{enumerate}
    \item\label{sim:setup-3} The errors are distributed according to a standard normal distribution with block-based contamination as presented in~\cite{zhang2013}. Here, in each resample, the normally distributed errors are first added to each of the entries of $\bb{X}$, and then a $2\times2$ block of submatrix is chosen randomly and all entries of that submatrix are substituted to $25$.
    \item\label{sim:setup-4} The errors are distributed according to a standard Cauchy distribution. This setup helps us to study the effect of heavy tailed errors in the robust estimation of SVD.
    \item\label{sim:setup-5} The errors are distributed according to a standard lognormal distribution, which is used to study the effect of an asymmetric error distribution with only positive support.
\end{enumerate}

Table~\ref{tab:summary} summarizes the comparative results of the usual SVD method, the three existing robust SVD algorithms (pcaSVD~\citep{pcaMethods} and two variants of RobRSVD~\citep{RobRSVD}) and the rSVDdpd method~\citep{roy2021new} for different choices of robustness parameter $\alpha$, based on the aforementioned performance measures, for different simulation setups~\ref{sim:setup-1}-\ref{sim:setup-5}. 


\begin{table}[tb]
    \centering
    \caption{Summary of performance measures (Total squared bias, Total MSE of singular values and Total dissimilarity, denoted by Diss, of left and right singular vectors) of different existing SVD and robust SVD algorithms}
    \label{tab:summary}
    \resizebox{\linewidth}{!}{
    \begin{tabular}{@{}crrrrrrrrrr@{}}
        \toprule
        \multirow{2}{*}{\begin{tabular}[c]{@{}l@{}}Simulation \\ Setup\end{tabular}} & \multirow{2}{*}{Measure} & \multicolumn{4}{c}{Existing methods for computing SVD} & \multicolumn{5}{c}{Choice of $\alpha$ in rSVDdpd} \\
         &  & \multicolumn{1}{c}{Usual SVD} & \multicolumn{1}{c}{pcaSVD} & \multicolumn{1}{c}{RobSVD} & \multicolumn{1}{c}{RobRSVD} & $\alpha = 0.1$ & $\alpha = 0.3$ & $\alpha = 0.5$ & $\alpha = 0.7$ & $\alpha = 1$ \\ \midrule
        \multirow{4}{*}{S1} & Sq. Bias &  7.957 & 15.066 & 8.808 & 1.684 & 7.952 & 7.94 & 7.925 & 7.894 & 7.764 \\  & Total MSE &  10.456 & 24.242 & 11.608 & 6.865 & 10.455 & 10.473 & 10.553 & 10.68 & 10.841 \\  & Diss (left) &  0.701 & 1.198 & 0.799 & 0.733 & 0.702 & 0.707 & 0.717 & 0.734 & 0.771 \\  & Diss (right) &  0.418 & 0.968 & 0.52 & 0.514 & 0.419 & 0.425 & 0.433 & 0.449 & 0.487 \\ \midrule 
        \multirow{4}{*}{S2a} & Sq. Bias &  519.35 & 350.402 & 523.859 & 17.672 & 294.047 & 18.138 & 10.877 & 9.433 & 8.779 \\  & Total MSE &  729.83 & 793.612 & 748.043 & 114.802 & 622.819 & 91.894 & 44.016 & 32.329 & 27.144 \\  & Diss (left) &  1.93 & 1.679 & 1.933 & 1.324 & 1.613 & 0.944 & 0.896 & 0.885 & 0.904 \\  & Diss (right) &  1.529 & 1.286 & 1.516 & 1.074 & 1.163 & 0.615 & 0.576 & 0.571 & 0.59 \\ \midrule
        \multirow{4}{*}{S2b} & Sq. Bias &  1397.36 & 1086.341 & 1427.869 & 60.112 & 1039.977 & 130.472 & 44.154 & 30.039 & 23.262 \\  & Total MSE &  1661.82 & 1706.458 & 1722.445 & 278.585 & 1529.374 & 494.34 & 237.68 & 167.244 & 132.886 \\  & Diss (left) &  2.155 & 1.965 & 2.148 & 1.672 & 2.013 & 1.283 & 1.113 & 1.084 & 1.078 \\  & Diss (right) &  1.702 & 1.5 & 1.665 & 1.371 & 1.5 & 0.899 & 0.756 & 0.733 & 0.731 \\ \midrule
        \multirow{4}{*}{S2c} & Sq. Bias &  2434.021 & 2110.99 & 2488.98 & 193.664 & 2094.11 & 662.604 & 298.008 & 196.514 & 141.064 \\  & Total MSE &  2712.69 & 2782.575 & 2803.595 & 617.094 & 2587.792 & 1453.404 & 938.705 & 739.152 & 612.155 \\  & Diss (left) &  2.206 & 2.11 & 2.203 & 1.932 & 2.17 & 1.666 & 1.444 & 1.363 & 1.313 \\  & Diss (right) &  1.73 & 1.633 & 1.688 & 1.598 & 1.651 & 1.227 & 1.042 & 0.965 & 0.93 \\ \midrule
        \multirow{4}{*}{S3} & Sq. Bias &  1677.949 & 1640.708 & 1679.881 & 1090.284 & 1355.176 & 114.714 & 28.128 & 21.594 & 18.048 \\  & Total MSE &  1686.361 & 1654.5 & 1688.708 & 1248.276 & 1634.409 & 458.468 & 156.986 & 113.33 & 96.675 \\  & Diss (left) &  2.052 & 2.003 & 1.941 & 2.162 & 2.007 & 1.208 & 1.024 & 1.012 & 1.002 \\  & Diss (right) &  1.924 & 1.836 & 1.832 & 2.175 & 1.844 & 0.895 & 0.692 & 0.678 & 0.667 \\ \midrule
        \multirow{4}{*}{S4} & Sq. Bias &  41825.788 & 14224.145 & 41779.81 & 540.799 & 469.522 & 265.135 & 198.809 & 169.082 & 140.645 \\  & Total MSE &  2171697.037 & 2163377.363 & 2171711.033 & 29149.731 & 1629.485 & 1197.881 & 859.039 & 807.869 & 842.77 \\  & Diss (left) &  2.089 & 1.989 & 2.095 & 1.707 & 1.97 & 1.877 & 1.838 & 1.809 & 1.786 \\  & Diss (right) &  1.603 & 1.489 & 1.602 & 1.303 & 1.496 & 1.403 & 1.367 & 1.355 & 1.324 \\ \midrule
        \multirow{4}{*}{S5} & Sq. Bias &  93.633 & 77.499 & 94.901 & 29.135 & 69.98 & 67.092 & 62.429 & 56.64 & 49.424 \\  & Total MSE &  146.187 & 108.514 & 149.13 & 49.046 & 85.102 & 83.272 & 78.83 & 72.46 & 65.682 \\  & Diss (left) &  2.02 & 2.034 & 1.969 & 1.967 & 1.968 & 1.954 & 1.943 & 1.939 & 1.934 \\  & Diss (right) &  1.785 & 1.81 & 1.734 & 1.824 & 1.754 & 1.744 & 1.734 & 1.731 & 1.726 \\
        \bottomrule
    \end{tabular}}
\end{table}

As shown in Table~\ref{tab:summary}, the usual SVD generally leads to a biased estimator of the singular values for Gaussian errors which is also supported in well-established theory~\citep{rudelson2010non}. For the simulation setup~\ref{sim:setup-1}, the pcaSVD algorithm is found to be the most biased, followed by RobSVD, both of which have more bias and MSE than the usual SVD algorithm. Compared to the usual SVD, rSVDdpd algorithm achieves lesser bias as the robustness parameter $\alpha$ increases, but at the cost of higher variance and MSE. RobRSVD achieves the minimum bias and MSE in this scenario, but it shows a higher variance and a higher dissimilarity score in singular vectors than the rSVDdpd algorithm.

Turning our attention to simulation setups~\ref{sim:setup-2a}, \ref{sim:setup-2b} and \ref{sim:setup-2c}, we see that the usual SVD and existing robust SVD algorithms pcaSVD and RobSVD do not yield very reliable estimates of the singular values. Although RobRSVD provides reasonable estimates, rSVDdpd achieves lower bias and MSE for some choices of $\alpha$. In the presence of random outlying observations, as in the case of simulation setups~\ref{sim:setup-2a}, \ref{sim:setup-2b} and \ref{sim:setup-2c}, both the bias and MSE for rSVDdpd show reductions as the robustness parameter $\alpha$ is increased from $0$ to $1$. The dissimilarities of singular vectors also tend to decrease with an increase in $\alpha$. 

For the block level contamination in simulation setup~\ref{sim:setup-3}, we find that rSVDdpd has much better performances than the other robust SVD algorithms for all performance metrics. With errors from a heavy-tailed distribution as considered in the simulation setup~\ref{sim:setup-4}, the results remain very similar. The rSVDdpd algorithm provides the least bias and MSE, and even with a small robustness parameter $\alpha = 0.1$, rSVDdpd outperforms the existing robust SVD algorithms under consideration. 

In simulation setup~\ref{sim:setup-5} with lognormally distributed errors having positive support, rSVDdpd outperforms the usual SVD, pcaSVD and RobSVD methods by showing a reduction in both bias and MSE. However, as in the simulation setup~\ref{sim:setup-1}, RobRSVD is again found to provide estimates with the least bias and MSE, but at a cost of higher variance and dissimilarity scores than rSVDdpd.

Although RobRSVD outputs better singular value estimates than rSVDdpd under normally and lognormally distributed errors, it does so at the cost of extremely high computational complexity. This is precisely due to the matrix inversion step to compute $(\mathcal{V}\tr \mathcal{W}^\ast \mathcal{V} + 2\Omega_{\bb{u}\mid \bb{v}})^{-1}$ (see Eq.~(9) of~\cite{zhang2013}). Since the best known matrix inversion algorithm, i.e., a variant of Coppersmith-Winograd algorithm~\citep{alman2021refined} achieves a computational complexity of $\mathcal{O}(n^{2.3728596})$ for inverting an $n \times n$ matrix, it follows that each iteration of the RobRSVD algorithm has time complexity $\mathcal{O}(n^{2.3728596} + p^{2.3728596})$. On the other hand, each iteration of rSVDdpd performs only a weighted average computation, which reduces its computational budget to $\mathcal{O}(n^2 + p^2)$. To demonstrate this, we consider $n \times p$ matrices with uniformly distributed entries for different choices of $n$ and fixed $p = 25$, and apply different methods of computing SVD on them. Table~\ref{tab:time-complexity} summarizes the time taken (in units of milliseconds) to obtain the first singular value from different algorithms for different choices of $n$, in a computer with Intel i5-8300H 2.30GHz processor with 8 GB of RAM. As seen from Table~\ref{tab:time-complexity}, the computational budget of rSVDdpd is similar to pcaSVD, which is lower by several orders of magnitude than RobSVD and RobRSVD. This extremely high computational cost of RobRSVD can be circumvented if the matrix $(\mathcal{V}\tr \mathcal{W}^\ast \mathcal{V} + 2\Omega_{\bb{u}\mid \bb{v}})$ becomes a diagonal matrix, which happens if the penalty parameter is taken as zero and RobRSVD is reduced to its non-regularized variant RobSVD. However, as Table~\ref{tab:summary} shows, the RobSVD algorithm without the regularization cannot provide a reliable robust estimate of singular values, even using Huber's robust loss function.

\begin{table}[tb]
    \centering
    \caption{Summary of average time taken (in milliseconds) to obtain the first singular value and vectors of an $n \times 25$ matrix with random entries from $U(0,1)$ via different SVD algorithms}
    \label{tab:time-complexity}
    \resizebox{\linewidth}{!}{
        \begin{tabular}{@{}crrrrrrr@{}}
            \toprule
            \multirow{2}{*}{Number of rows ($n$)} & \multicolumn{4}{c}{Existing methods for computing SVD} & \multicolumn{3}{c}{Choice of $\alpha$ in rSVDdpd} \\
            & \multicolumn{1}{c}{Usual SVD} & \multicolumn{1}{c}{pcaSVD} & \multicolumn{1}{c}{RobSVD} & \multicolumn{1}{c}{RobRSVD} & $\alpha = 0.1$ & $\alpha = 0.5$ & $\alpha = 1$ \\ \midrule
            5 & 0.014 & 2.209 & 3.098 & 73.666 & 0.545 & 0.841 & 1.417 \\
            10 & 0.011 & 4.388 & 3.966 & 99.077 & 0.975 & 1.508 & 2.100 \\
            25 & 0.008 & 4.965 & 3.989 & 81.045 & 3.861 & 6.416 & 11.394 \\
            50 & 0.013 & 8.519 & 7.824 & 149.4 & 4.396 & 7.104 & 12.526 \\
            100 & 0.017 & 15.696 & 34.649 & 826.44 & 5.136 & 8.683 & 14.362 \\
            250 & 0.026 & 32.066 & 402.125 & 7839.942 & 7.756 & 13.252 & 22.379 \\
            500 & 0.041 & 35.828 & 2948.001 & 54209.15 & 13.622 & 21.413 & 32.404 \\
            750 & 0.058 & 58.697 & 10363.441 & 210494.564 & 24.67 & 36.563 & 55.422 \\
            1000 & 0.072 & 69.76 & 26282.893 & 531362.110 & 27.727 & 40.309 & 62.234 \\
            \bottomrule
        \end{tabular}
    }
\end{table}

\section{Conclusion}\label{sec:conclusion}

As depicted in Section~\ref{sec:intro}, a plethora of algorithms from an extensive range of disciplines use singular value decomposition as a core component of the methods. However, the increasing prevalence of big data has made it challenging to ensure the accuracy and reliability of the data. The input data for these algorithms are prone to contamination by noise and outliers, leading to inaccurate results when using standard SVD. To address this issue, several robust SVD methods have been proposed (see Section~\ref{sec:related-works}), but most of them are not scalable to large matrices encountered in real-life applications. The lack of theoretical guarantees of these algorithms has limited their widespread adoption and hindered their application in critical domains. While \cite{roy2021new} demonstrate an application of the ``rSVDdpd'' algorithm to solve a real-life problem, in this paper, we provide the theoretical justification for its reliability. The simulation results further validate the superiority of rSVDdpd compared to the existing algorithms. However, more investigation is needed to develop asymptotic distributions of the estimated singular values and vectors, which can provide confidence interval estimates.

We believe that the ``rSVDdpd'' algorithm has potential applications beyond video surveillance background modelling by using it as a replacement of the standard SVD in various algorithms to handle data contamination. As an example, we can use rSVDdpd for modelling genetic data, performing community detection in networks, estimating latent semantic representation of text documents from term-document matrices, etc. We hope to explore these applications in future.

\appendix

\section{A Brief Review of Minimum Density Power Divergence Estimator}\label{appendix:dpd}

    \cite{basu-dpd} introduced the density power divergence as a measure of discrepancy between two probability density functions, which being an M-estimator, as well as a minimum distance estimator, enjoys various theoretical properties. The density power divergence between the densities $g$ and $f$ is defined as
    \begin{equation*}
        d_{\alpha}(g, f) = \begin{cases}
        \displaystyle\int \left\{ f^{1+\alpha} - \left(1 + \dfrac{1}{\alpha} \right)f^\alpha g + \dfrac{1}{\alpha} g^{1+\alpha} \right\} & \alpha > 0\\
        \displaystyle\int g \ln\left( \dfrac{g}{f} \right) & \alpha = 0\\
        \end{cases}.    
    \end{equation*}
    \noindent Here $\ln(\cdot)$ denotes the natural logarithm. The control parameter $\alpha$ provides a smooth bridge between robustness and efficiency. 
    
    In case of independent and identically distributed observations, $Y_1, Y_2, \dots Y_n$, with true distribution function $G$ and corresponding density $g$, we model this unknown density by a parametric family of densities $\mathcal{F}_{\bb{\theta}} = \left\{ f_{\bb{\theta}} : \bb{\theta} \in \bb{\Theta} \right\}$. The estimator of $\theta$ is then obtained as
    \begin{equation*}
        \widehat{\bb{\theta}} = \arg\min_{\bb{\theta} \in \bb{\Theta} } d_{\alpha}(dG_n, f_{\bb{\theta}}) ,   
    \end{equation*}
    \noindent where $G_n$ is the empirical distribution function. This can be shown to be equivalent to
    \begin{equation*}
        \widehat{\bb{\theta}} = \arg\min_{\bb{\theta} \in \bb{\Theta} } \left[  \int f_{\bb{\theta}}^{1 + \alpha} - \left( 1 + \dfrac{1}{\alpha} \right)\dfrac{1}{n}\sum_{i=1}^{n} f_{\bb{\theta}}(Y_i)^{\alpha} \right].
    \end{equation*}
    \noindent Later, \cite{ghosh2013} extended this work by allowing independent but not identically distributed data. In this case, the observed data $Y_i \sim g_i$, where each $g_i$ is an unknown density. Each of the true density $g_i$ is modeled by a corresponding parametric family of densities $\mathcal{F}_{i, \theta} = \left\{ f_{i, \theta} : \theta \in \bb{\Theta} \right\}$ for all $i = 1, 2, \dots n$. Finally, the proposed MDPD estimator is obtained as
    \begin{equation*}
        \widehat{\bb{\theta}} = \arg\min_{\bb{\theta} \in \bb{\Theta} } \dfrac{1}{n} \sum_{i=1}^{n} \left[ \int f_{i, \bb{\theta}}^{1 + \alpha} - \left( 1 + \dfrac{1}{\alpha} \right)f_{i, \bb{\theta}}(Y_i)^{\alpha} \right].
    \end{equation*}
    \noindent Various nice theoretical properties like consistency and asymptotic normality of the above MDPD estimator have been proven by~\cite{ghosh2013}.
    
\section{Proofs of the Results}\label{appendix:proofs}

\subsection{Proof of Theorem~\ref{thm:convergence}}\label{appendix:convergence-proof}

Note that, each $V_{ij,\alpha}^{(1)}(\bb{\theta})$ is bounded below by the finite quantity $\epsilon^{-\alpha}(M_f - (1+1/\alpha)f^\alpha(0))$, hence the same lower bound also applies for $H_\alpha^{(1)}(\bb{\theta})$. Therefore, there exists at least one local minimum of $H_\alpha^{(1)}(\bb{\theta})$.

We first show that iterating equations~\eqref{eqn:iteration-rule-1}-\eqref{eqn:iteration-rule-3} reduces the value of the objective function $H_\alpha^{(1)}(\bb{\theta})$. We shall show this only for Eq.~\eqref{eqn:iteration-rule-1}, rest can be shown similarly. Let, $e_{ij}^{(t+1/2)} = X_{ij} - \lambda^{(t+1)}u_i^{(t+1)}v_j^{(t)}$. Then,
\begin{equation*}
    e_{ij}^{(t)} = e_{ij}^{(t+1/2)} + v_j^{(t)} \dfrac{\sum_{k} v_k^{(t)} e_{ik}^{(t)} \psi(e_{ik}^{(t)}/\sigma^{(t)}) }{\sum_k (v_k^{(t)})^2 \psi(e_{ik}^{(t)}/\sigma^{(t)}) }.
\end{equation*}
\noindent Let's call the second term in the above sum as $v_j^{(t)} e_{i}^\ast$. An application of Cauchy-Schwartz inequality along with Assumption~\ref{assum:psi-bound} shows that $\vert e_i^\ast\vert \leq K$ for some constant $K$. Together with $\vert v_j^{(t)}\vert \leq 1$, it ensures that there exists $K_1, K_2 > 0$ such that
\begin{equation*}
    \vert e_{ij}^{(t)} \vert \leq K_1 + K_2 t.
\end{equation*}
\noindent In view of the definition of $H_\alpha^{(1)}(\bb{\theta})$, it is now enough to show that
\begin{equation*}
    \sum_{i,j} \left[ f^\alpha\left( \dfrac{\vert e_{ij}^{(t+1/2)} \vert }{\sigma^{(t)}} \right) - f^\alpha\left( \dfrac{\vert e_{ij}^{(t)} \vert }{\sigma^{(t)}} \right)  \right] \geq 0.
\end{equation*}
\noindent An application of Taylor's theorem yields
\begin{equation*}
    f^\alpha\left( \dfrac{\vert e_{ij}^{(t+1/2)} \vert }{\sigma^{(t)}} \right)
    - f^\alpha\left( \dfrac{\vert e_{ij}^{(t)} \vert }{\sigma^{(t)}} \right) 
    = -\dfrac{\alpha}{\sigma^{(t)}} \psi\left( \dfrac{\vert e_{ij}^{(t)} \vert }{\sigma^{(t)}} \right)(\vert e_{ij}^{(t)} - v_j^{(t)} e_{i}^\ast \vert - \vert e_{ij}^{(t)} \vert) 
    - \dfrac{\alpha}{2 (\sigma^2)^{(t)} } \psi'(c) (\vert e_{ij}^{(t)} - v_j^{(t)} e_{i}^\ast \vert - \vert e_{ij}^{(t)} \vert)^2,
\end{equation*}
\noindent where $c$ is some value between $e_{ij}^{(t+1/2)}$ and $e_{ij}^{(t)}$. Because of Assumption~\ref{assum:f-diff-cond}, we have $\psi'(x) > 0$ and hence the second term is nonnegative. For the first term, consider the inequality
\begin{equation*}
    \sum_{j} \psi\left( \dfrac{\vert e_{ij}^{(t)} \vert }{\sigma^{(t)}} \right)(\vert e_{ij}^{(t)} \vert - \vert e_{ij}^{(t)} - v_j^{(t)} e_{i}^\ast \vert ) 
    \geq \sum_{j} \psi\left( \dfrac{\vert e_{ij}^{(t)} \vert }{\sigma^{(t)}} \right) \dfrac{(e_{ij}^{(t)})^2 - (e_{ij}^{(t)} - v_j^{(t)} e_{i}^\ast )^2 }{2K_1 + K_2(2t+1)},
\end{equation*}
\noindent since, $\vert e_{ij}^{(t)} \vert + \vert e_{ij}^{(t)} - v_j^{(t)} e_{i}^\ast \vert \leq 2K_1 + K_2(2t+1)$. This lower bound is nonnegative since by the structure of $e_i^\ast$, it minimizes the weighted squared error
\begin{equation*}
    \sum_{j} \psi\left( \dfrac{\vert e_{ij}^{(t)} \vert }{\sigma^{(t)}} \right) (e_{ij}^{(t)} - v_j^{(t)} a )^2,
\end{equation*}
\noindent over choice of all possible $a \in \R$. Adding these quantities for all $i$ and putting it back to Taylor's series shows that each iteration decreases the value of the objective function $H_\alpha^{(1)}$. Now, the sequence $\{ H_\alpha^{(1)}(\bb{\theta}^{(t)} \}_{t=0}^\infty$ becomes a decreasing sequence of real numbers bounded below, and hence has a convergent subsequence. Then, the facts that $H_\alpha^{(1)}$ is a continuous function of $\bb{\theta}$ due to continuity of $f$ and that $\bb{\Theta}$ can effectively be restricted to a compact set $[0, \Vert \bb{X}\Vert_F] \times S_n^+ \times S_p \times [\epsilon, \Vert \bb{X}\Vert_F]$ imply that $\bb{\theta}^{(t)}$ converges to some $\bb{\theta}^\ast$. Finally, since $\bb{\theta}^\ast$ satisfies the iterating equations~\eqref{eqn:iteration-rule-1}-\eqref{eqn:iteration-rule-3}, it in turn, satisfies the estimating equations, i.e., the gradient of $H_\alpha^{(1)}$ at $\bb{\theta}^\ast$ is zero. This implies that $\bb{\theta}^\ast$ is a local minimum of the same.

\subsection{Proof of Theorem~\ref{thm:def-match}}

Relation~\eqref{eqn:true-norm} is verified by the implications that $u_i^\ast$ and $v_j^\ast$s belong to the respective Stiefel manifold. To verify~\eqref{eqn:true-ai}, note that with $v_j^g = v_j^\ast$ and $\sigma^g = (\sigma^\ast)$, the quantity in~\eqref{eqn:true-ai} is same as minimizing
\begin{equation*}
    \int V_f(x; a, v_j^\ast, \sigma^\ast) + \dfrac{1}{\alpha} \int g_{ij}^\alpha,
\end{equation*}
\noindent since the last term is independent of the minimization over $a$. But this is the density power divergence (DPD) between the density $f$ and true density $g_{ij}$. From Theorem 2.1 of~\cite{basu-dpd}, it follows that this divergence is minimized if and only if two densities match, i.e., $a = \lambda^\ast u_i^\ast$. By exactly similar logic and interchanging the roles of $u_i$ and $v_j$,~\eqref{eqn:true-bj} and~\eqref{eqn:true-sigma} can also be verified. This proves that $\bb{\theta}^\ast = ( \lambda^\ast, \{ u_i^\ast \}_{i=1}^n, \{ v_j^\ast \}_{j = 1}^p, (\sigma^\ast)^2 )$ is a ``best'' fitting parameter for the given setup.

In order to prove uniqueness, suppose $\widetilde{\bb{\theta}} = ( \widetilde{\lambda}, \{ \tilde{u}_i \}_{i=1}^n, \{ \tilde{v}_j \}_{j = 1}^p, \tilde{\sigma}^2 )$ be another ``best'' fitting parameter. Then again, the DPD with $v_j$ and $\sigma$ substituted for $v_j^\ast$ and $\sigma^\ast$ respectively, is minimized at $a = \widetilde{\lambda} \tilde{u}_i$ independently of the choice of $j$. However, this divergence can be made equal to its minimum value $0$ if and only if $\tilde{\sigma}^2 = (\sigma^\ast)^2$, and
\begin{equation}
    \lambda u_i v_j = \lambda^\ast u_i^\ast v_j^\ast, \ i = 1, \dots n; j = 1, \dots p,\label{eqn:unique-lambda}
\end{equation}
\noindent which follows from Theorem 2.1 of~\cite{basu-dpd}. Since, both $\widetilde{\bb{\theta}}$ and $\bb{\theta}^\ast$ are ``best'' fitting parameters, they must satisfy~\eqref{eqn:true-norm}. Hence, $(\widetilde{\lambda} \tilde{u}_i)^2 = \sum_j (\widetilde{\lambda} \tilde{u}_i \tilde{v}_j)^2 = \sum_j (\lambda^\ast u_i^\ast v_j^\ast)^2 = (\lambda^\ast u_i^\ast)^2$. Taking sum over the row index $i$ now gives $\widetilde{\lambda}^2 = (\lambda^\ast)^2$. Since both $\widetilde{\lambda}, \lambda^\ast \geq 0$, it follows that $\widetilde{\lambda} = \lambda^\ast$, and consequently, $\vert \tilde{u}_i \vert = \vert u_i^\ast \vert$ and $\vert \tilde{v}_j \vert = \vert v_j^\ast \vert$.

Now suppose $\tilde{v}_j = (-v_j^\ast)$ for some $j$. Along with~\eqref{eqn:unique-lambda}, it means that $\tilde{u}_i = (-u_i^\ast)$ for all $i = 1, 2, \dots n$. This leads to a contradiction since both $\{ \tilde{u}_i \}$ and $\{ u_i^\ast \}$ cannot be in $S_n^+$.

\subsection{Proof of Theorem~\ref{thm:scale-equivariance}}

It is obvious that $\widetilde{\bb{\theta}}^g$ satisfies~\eqref{eqn:true-norm} as $\bb{\theta}^g$ is given to be a ``best'' fitting parameter. Considering the matrix $\bb{Y} = c\bb{X}$, let us denote the true density of $Y_{ij}$ as $g_{ij}^Y(\cdot)$, as opposed to $g_{ij}(\cdot)$ denoting the true density of $X_{ij}$. A change of variable formula yields that $g_{ij}^Y(y) = g_{ij}(y / c)$. Hence, from the substitution principle of integration, it follows that
\begin{equation*}
    \int V_f\left(y; a, v_j^g, c^2(\sigma^g)^2\right)g_{ij}^Y(y) dy
    = c^\alpha \int V_f\left(z; a/c, v_j^g, (\sigma^g)^2\right) g_{ij}(z)dz.
\end{equation*}
\noindent Since the right-hand side is minimized at $a / c = \lambda^g u_i^g$, the left-hand side is minimized at $a = c \lambda^g u_i^g$. This verifies~\eqref{eqn:true-ai}. Similar to this,~\eqref{eqn:true-bj} can also be established by interchanging the role of $\bb{a}$ and $\bb{b}$ above. For the parameter $\sigma$, again a substitution principle applies, and we obtain
\begin{equation*}
    \int V_f(y; c\lambda^ga_i^g, b_j^g, \sigma^2)g_{ij}^Y(y) dy 
    = c^\alpha \int V_f(z; \lambda a_i^g, b_j^g, \sigma^2 / c^2)g_{ij}(z) dz.
\end{equation*}
Again by the hypothesis that $\bb{\theta}^g$ is a ``best'' fitting parameter, the latter is minimized when $\sigma/c = \sigma^g$, hence the former is minimized at $\sigma^2 = c^2 (\sigma^g)^2$. This verifies~\eqref{eqn:true-sigma}.

\subsection{Proof of Theorem~\ref{thm:permute-equivariance}}

Let, $\bb{Y} = \bb{PXQ}\tr$. Then,~\eqref{eqn:true-norm} is satisfied for the new setup as $\sum_{i=1}^{n} (u_{\pi_R(i)}^g)^2 = \sum_{i=1}^{n} (u_{i}^g)^2 = 1$ and similarly $\sum_{j=1}^{p} (v_{\pi_C(j)}^g)^2 = \sum_{j=1}^{p} (v_{j}^g)^2 = 1$. To see that~\eqref{eqn:true-ai} hold for the new setup with $\widetilde{\theta}^g$, note that for every $j = 1, 2, \dots p$, the minimizer of the integral in~\eqref{eqn:true-ai} is $u_i^g$, independent of the choice of $j$. Now, considering~\eqref{eqn:true-ai} for $\pi_R^{-1}(i)$ and $\pi_C^{-1}(j)$ instead of $i$ and $j$, we obtain
\begin{equation}
    \lambda^g u_{\pi_R^{-1}(i)}^g = \argmin_{a} \int V_f(x; a, v_{\pi_C^{-1}(j)}^g, (\sigma^g)^2 ) 
    g_{\pi_R^{-1}(i), \pi_C^{-1}(j)}(x) dx.
    \label{eqn:theorem-2-proof-1}
\end{equation}
\noindent However, $u^{g}_{\pi_R^{-1}(i)}$ is the $i$-th entry of sequence $\{ u_{\pi_R(i)}^g : i = 1, 2, \dots n \}$, i.e. if we consider a vector $\bb{u}^g$ with its entries $u_i^g$, then $u^{g}_{\pi_R^{-1}(i)}$ is the $i$-th entry of $\bb{Pu}$. Similarly, $v^{g}_{\pi_C^{-1}(j)}$ is the $j$-th entry of $\bb{Qv}$. And finally, the elements of the new matrix are $Y_{ij} = X_{\pi_R(i), \pi_C(j)}$, thus the density for the element $Y_{ij}$ is $g^Y_{ij}(y) = g_{\pi^{-1}_R(i), \pi^{-1}_C(j)}(y)$ which can be verified by a change of variable formula. Combining these,~\eqref{eqn:theorem-2-proof-1} can be reformulated as
\begin{equation*}
    \lambda^g (\bb{Pu}^g)_{i}
    = \argmin_{a} \int V_f(x; a, (\bb{Qv}^g)_{j}, (\sigma^g)^2) g^Y_{ij}(x) dx.
\end{equation*}
\noindent This shows that~\eqref{eqn:true-ai} holds for new matrix $\bb{PXQ}\tr$ with the given best fitting parameter $\widetilde{\bb{\theta}}^g$. The relation~\eqref{eqn:true-bj} holds by imitating the same proof, except interchanging the role of $\bb{a}$ and $\bb{b}$. Finally, relation~\eqref{eqn:true-sigma} for the permuted matrix follows from noting that
\begin{equation*}
    \int V_f(x; \lambda^g u_{\pi_R^{-1}(i)}^g, v_{\pi_C^{-1}(j)}^g, \sigma^2) g_{\pi_R^{-1}(i), \pi_C^{-1}(j)}(x) dx\\
    = \int V_f\left(x; \lambda^g (\bb{Pu}^g)_{i}, (\bb{Qv}^g)_{j}, \sigma^2 \right) g^Y_{ij}(x) dx.
\end{equation*}

\subsection{Proof of Theorem~\ref{thm:scale-equivariance-estimate}}\label{appendix:scale-equiv-estimate-proof}

Let us denote $\bb{\theta}^{(t)}$ denote the estimate at $t$-th iteration for data matrix $\bb{X}$ and let $\widetilde{\bb{\theta}}^{(t)}$ denote the same for the matrix $c\bb{X}$. Clearly, it is then enough to show that for all $t = 1, 2, \dots$,
\begin{equation*}
    \widetilde{\lambda}^{(t)} = c \lambda^{(t)}, \
    (\widetilde{\sigma}^{(t)})^2 = c^2 (\sigma^{(t)})^2, 
    \widetilde{u}_i^{(t)} = u_i^{(t)}, \ \widetilde{v}_j^{(t)} = v_j^{(t)}.
\end{equation*}
\noindent We will show this by using the principle of mathematical induction. For $t = 0$, the claim is validated by the equivariance of the initial estimate. To show the inductive step, we first consider Eq.~\eqref{eqn:iteration-rule-1}. Note that,
\begin{equation*}
    \widetilde{e}_{ij}^{(t)} = c X_{ij} -  \widetilde{\lambda}^{(t)}\widetilde{u}_i^{(t)}\widetilde{v}_j^{(t)} = c e_{ij}^{(t)},
\end{equation*}
\noindent by induction hypothesis. Therefore,
\begin{equation*}
    \widetilde{\lambda}^{(t+1)} \widetilde{u}_i^{{(t+1)}}
    = \dfrac{\sum_j c\widetilde{v}_j^{(t)} X_{ij} \psi(\widetilde{e}_{ij}^{(t)}/\widetilde{\sigma}^{(t)}) }{\sum_j (\widetilde{v}_j^{(t)})^2 \psi(\widetilde{e}_{ij}^{(t)}/\widetilde{\sigma}^{(t)}) } \\
    = c \dfrac{\sum_j v_j^{(t)} X_{ij} \psi(c e_{ij}^{(t)}/c{\sigma}^{(t)}) }{\sum_j ({v}_j^{(t)})^2 \psi(c{e}_{ij}^{(t)}/c{\sigma}^{(t)}) }                                                           \\
    = c {\lambda}^{(t+1)} {u}_i^{{(t+1)}}.
\end{equation*}
\noindent Performing the same steps with Eq.~\eqref{eqn:iteration-rule-2} and~\eqref{eqn:iteration-rule-3} ensure that
\begin{equation*}
    \widetilde{\lambda}^{(t+1)} \widetilde{v}_j^{{(t+1)}} = c{\lambda}^{(t+1)} {v}_j^{{(t+1)}}, \ \widetilde{\sigma}^{(t+1)} = c\sigma^{(t+1)}.
\end{equation*}
\noindent Finally, since the estimates of the singular vectors are normalized and restricted to be in the parameter space $\bb{\Theta} = [0, \infty) \times S_n^+ \times S_p^+ \times [0, \infty)$, the inductive step follows from a normalization step.

\subsection{Proof of Theorem~\ref{thm:permute-equivariance-estimate}}

This proof is very similar to the proof of Theorem~\ref{thm:scale-equivariance-estimate}. We shall again denote $\bb{\theta}^{(t)}$ as the estimate at the $t$-th iteration for the data matrix $\bb{X}$ and $\widetilde{\bb{\theta}}^{(t)}$ as the estimate at the $t$-th iteration for the data matrix $\bb{PXQ}\tr$. Again, it is enough to show that 
\begin{align*}
    \widetilde{\lambda}^{(t)} = \lambda^{(t)}, \ 
    (\widetilde{\sigma}^{(t)})^2 = (\sigma^{(t)})^2, \
    \widetilde{u}_i^{(t)} = u_{\pi_R(i)}^{(t)}, \ \widetilde{v}_j^{(t)} = v_{\pi_C(j)}^{(t)}, \quad
    i = 1, \dots n; \ j = 1, \dots p; \ t = 1, 2, \ldots,
\end{align*}
\noindent which we shall show using the principle of mathematical induction. The initial case $t = 0$ follows from the equivariance of the initial estimate. To show the inductive step, note that
\begin{equation*}
    \widetilde{e}_{ij}^{(t)} = X_{\pi_R(i), \pi_C(j)} -  \widetilde{\lambda}^{(t)}\widetilde{u}_{\pi_R(i)}^{(t)}\widetilde{v}_{\pi_C(j)}^{(t)} = e_{\pi_R(i), \pi_C(j)}^{(t)}.
\end{equation*}
\noindent Now considering Eq.~\eqref{eqn:iteration-rule-1}, we get that 
\begin{align*}
    \widetilde{\lambda}^{(t+1)} \widetilde{u}_i^{{(t+1)}}
    & = \dfrac{\sum_j \widetilde{v}_j^{(t)} \widetilde{X}_{ij} \psi(\widetilde{e}_{ij}^{(t)}/\widetilde{\sigma}^{(t)}) }{\sum_j (\widetilde{v}_j^{(t)})^2 \psi(\widetilde{e}_{ij}^{(t)}/\widetilde{\sigma}^{(t)}) } \\
    & = \dfrac{\sum_j v_{\pi_C(j)}^{(t)} X_{\pi_R(i), \pi_C(j)} \psi(c e_{\pi_R(i), \pi_C(j)}^{(t)}/{\sigma}^{(t)}) }{\sum_j ({v}_{\pi_C(j)}^{(t)})^2 \psi({e}_{\pi_R(i), \pi_C(j)}^{(t)}/{\sigma}^{(t)}) } \\
    & = \dfrac{\sum_{j'} v_{j'}^{(t)} X_{\pi_R(i), j'} \psi(c e_{\pi_R(i), j'}^{(t)}/{\sigma}^{(t)}) }{\sum_{j'} ({v}_{j'}^{(t)})^2 \psi({e}_{\pi_R(i), j'}^{(t)}/{\sigma}^{(t)}) } , \text{ calling the index } \pi_C(j) \text{ as } j'\\
    & = {\lambda}^{(t+1)} {u}_{\pi_R(i)}^{{(t+1)}}.
\end{align*}
\noindent We can perform the same steps with Eq.~\eqref{eqn:iteration-rule-2} and~\eqref{eqn:iteration-rule-3} as well, which completes the inductive step.

\subsection{Proof of Theorem~\ref{thm:consistency}}\label{appendix:consistency-proof}

First, we observe that the stereographic transformation mentioned in the discussion prior to Theorem~\ref{thm:consistency} can be employed and would remain valid because of Assumptions~\ref{assum:open-rect} and~\ref{assum:stereographic}. Now, to prove the consistency, we shall take a route similar to the one taken by \cite{ghosh2013} as in the case of MDPDE for INH setup. Instead of showing that the rSVDdpd estimator i.e., $\bb{\theta}^\ast_{n,p}$ is consistent for $\bb{\theta}^g_{n,p}$, we shall show instead that $\bb{\eta}^\ast_{n,p}$ is consistent for $\bb{\eta}^g_{n,p}$. Let us denote $\mathcal{H}_{n,p}(\bb{\eta})$ to indicate the $H$-function as in~\eqref{eqn:H-general} evaluated at $\bb{\theta} = \mathcal{T}^{-1}(\bb{\eta})$, for fixed $n$ and $p$ with $V_f$ substituted by $V_\phi$ given in~\eqref{eqn:V-normal-function}. To prove that $\bb{\eta}^\ast_{n,p}$ is consistent for $\bb{\eta}^g_{n,p}$, we shall show that for any sufficiently small $r > 0$, $\mathcal{H}_{n,p}(\bb{\eta}_{n,p}) > \mathcal{H}_{n,p}(\bb{\eta}^g_{n,p})$ for sufficiently large $n$ and $p$, for any $\bb{\eta}_{n,p}$ with $\Vert \bb{\eta}_{n,p} - \bb{\eta}^g_{n,p} \Vert_2 = r$. This means that the value of $\mathcal{H}_{n,p}$ at the surface of the ball of radius $r$ centered at $\bb{\eta}^g_{n,p}$ would be higher than its value at $\bb{\eta}^g_{n,p}$, and hence by the smoothness of $\mathcal{H}_{n,p}$, it is ensured that there will be a local minimum strictly inside that ball. Proceeding as in~\cite{ghosh2013}, we start with the Taylor series expansion of $\mathcal{H}_{n,p}(\bb{\eta}_{n,p})$ about $\bb{\eta}^g_{n,p}$, for any fixed $n$ and $p$. For notational convenience, we suppress the subscripts $n$ and $p$ from $\bb{\eta}$ and $\bb{\eta}^g$ which should be obvious from the context. We also use the symbol $\partial_{x_{i_1},\dots x_{i_k}}\Hcal_{n,p}$ to denote the $k$-th order partial derivative of $\Hcal_{n,p}$ in the direction of the variables $x_{i_1}, \dots x_{i_k}$ respectively, at the true parameter $\bb{\eta}^g$.
\begin{align*}
    & \mathcal{H}_{n, p}(\bb{\eta}) - \mathcal{H}_{n, p}(\bb{\eta}^g) \\
    = & \Hdiff{\lambda}(\lambda - \lambda^g) + \sum_{i=1}^{(n-1)} \Hdiff{\alpha_i}(\alpha_i - \alpha_i^g)  + \sum_{j=1}^{(p-1)}\Hdiff{\beta_j}(\beta_j - \beta_j^g) + \Hdiff{\sigma^2}(\sigma^2 - (\sigma^g)^2)    \\
    + & \dfrac{1}{2}\sum_{k_1, k_2} \Hdifftwo{\eta_{k_1},\eta_{k_2}} (\eta_{k_1} - \eta_{k_1}^g)(\eta_{k_2} - \eta_{k_2}^g) 
    + \frac{1}{6} \sum_{k_1,k_2,k_3} \partial^3_{\eta_{k_1},\eta_{k_2}, \eta_{k_3}}\Hcal_{n,p} (\eta_{k_1} - \eta_{k_1}^g) (\eta_{k_2} - \eta_{k_2}^g) (\eta_{k_3} - \eta_{k_3}^g)  \\
    = & S_{1,1} + S_{1,2} + S_{1, 3} + S_{1, 4} + \frac{1}{2}S_2 + \frac{1}{6} S_3,
\end{align*}
\noindent where the quantities $S_{1, 1}, S_{1, 2}, S_{1, 3}, S_{1,4}, S_2$ and $S_3$ respectively denote the summands they are replacing. Here, $\eta_k$ denotes the $k$-th coordinate of the vector $\bb{\eta}_{n,p}$. Also, $\alpha_i$'s and $\beta_j$'s are the natural parametric representation of the elements of left ($u_{i,n}$) and right singular vectors ($v_{j,p}$) respectively, where the dimension subscripts ($n$ and $p$) have been suppressed for notational convenience as indicated before.

Clearly, the smoothness of $\mathcal{H}_{n,p}$ along with Assumption~\ref{assum:best-fitting-distn} on the normality of the errors, indicates that $\int \mathcal{H}_{n,p} g_{ij}(x)dx$ can be differentiated thrice with respect to $\bb{\eta}_{n,p}$, and the derivative can be taken under the integral sign. Hence, we have
\begin{equation}
    \E\left[ \Hdiff{\eta_k} \right] = \partial_{\eta_k} \E\left(\Hcal_{n,p}\right) = 0,
    \label{eqn:proof-theorem5-1}
\end{equation}
\noindent since the population version of the objective function $\E\mathcal{H}_{n,p}$ is minimized at the true parameter $\bb{\eta}^g$. Thus, by a generalized version of Khinchin's Weak Law of Large numbers, it follows that as $n$ and $p$ both increase to infinity, each of the first order partial derivatives goes in probability to $0$. However, the problem arises as there are potentially infinitely many terms (as the parameter space increases in dimension). This jeopardizes any approach to naturally extending the proof of Theorem 3.1 of~\cite{ghosh2013}.

Before proceeding further, we note that since $\sum_{k=1}^n (u_k^g)^2 = 1$, its derivative yields $\sum_{k=1}^n a_u^g \partial_{\alpha_i}u_k = 0$ for any $i = 1, \dots (n-1)$. Similarly, $\sum_{l=1}^p b_l^g \partial_{\beta_j}v_l = 0$ for all $j = 1, \dots (p-1)$. Also, for notational convenience, we denote $w_{ij} = e^{-\alpha Z_{ij}^2/2}$.

Let us consider each of the sums $S_{1, 1}, S_{1, 2}, S_{1, 3}$ and $S_{1, 4}$ pertaining to the first order derivative separately. Since $\Hdiff{\lambda}$ and $\Hdiff{\sigma^2}$ both converges in probability to $0$, hence for sufficiently large $n$ and $p$, we have $\vert S_{1, 1} \vert < r^3$  and $\vert S_{1, 4} \vert < r^3$ with probability tending to $1$. Now, to deal with an increasing number of summands in $S_{1, 2}$ or $S_{1, 3}$ we apply Chebyshev's inequality after bounding its expectation and variance separately. By chain rule of differentiation,
\begin{equation}
    s_{n,p}
    = \sum_{i=1}^{(n-1)} \Hdiff{\alpha_i}
    = \sum_{k=1}^n \partial_{u_k}\Hcal_{n,p} \sum_{i=1}^{(n-1)} \partial_{\alpha_i}u_k,
    \label{eqn:proof-theorem5-2}
\end{equation}
\noindent where $\partial_{\alpha_i}u_k$ denotes the partial derivative of the entry of the left singular vector $u_k$ with respect to the stereographic projection variables $\alpha_i$ at $\bb{\eta}^g$. As in the case of~\eqref{eqn:proof-theorem5-1}, one can verify that $\E\left[ \Hdiff{u_k} \right] = 0$ for all $k = 1, 2, \dots n$, and, therefore,~\eqref{eqn:proof-theorem5-2} implies that $\E(s_{n,p}) = 0$. Turning to its variance, it follows that
\begin{equation*}
    \var(s_{n,p})                                                                                                                 = \sum_{k=1}^n \left( \sum_{i=1}^{(n-1)} \partial_{\alpha_i}u_k \right)^2 \var\left( \Hdiff{u_k} \right)  = \sum_{k=1}^n \left( \sum_{i=1}^{(n-1)} \partial_{\alpha_i}u_k \right)^2 \frac{(\alpha + 1)^2 (\lambda^g)^2 }{(2\pi)^{\alpha} \sigma^{2(\alpha + 1)} n^2p^2} B_1,
\end{equation*}
\noindent where $B_1 = \E(Z^2_{ij}w^2_{ij})$. Here we use the fact that for $k \neq l$, $\cov\left( \Hdiff{u_k}, \Hdiff{u_l} \right) = 0$, which follows by noting that the part of $\mathcal{H}_{n,p}$ dependent on $u_k$ would consist of only the $k$-th row of the data matrix $\bb{X}$, which are assumed to be independently distributed in the current setup. It also follows from Cauchy-Schwartz inequality that $S_u = \sum_{i=1}^{(n-1)} u_i^g \leq \sqrt{n-1}$, hence the sum
\begin{equation*}
    \sum_{k=1}^n \left( \sum_{i=1}^{(n-1)} \partial_{\alpha_i}u_k \right)^2 
    = (1 - u_n^g)^2 S_u^2 + \sum_{k=1}^n \left( (1 - u_n^g) - u_k^g S_u \right)^2 
    = (1 - u_n^g)^2 S_u^2 + n(1-u_n^g)^2 - S_u^2,
\end{equation*}
\noindent is bounded by $11n$ in magnitude. Therefore, for sufficiently large $n$ and $p$,
\begin{equation*}
    \var(s_{n,p}) = \mathcal{O}\left( (\sigma^g)^{-(2\alpha+2)} / np^2 \right).
\end{equation*}
\noindent Since we have $\sigma^g \asymp (np)^{-1/4}$ and $\alpha \leq 1$, it follows that $\var(s_{n,p}) \rightarrow 0$ as $n$ and $p$ tends to infinity. Therefore, we have $\vert \sum_{i=1}^{(n-1)} \Hdiff{\alpha_i} \vert \rightarrow 0$, with probability tending to one. Along with $\vert \alpha_i - \alpha_i^g\vert < r$, we have $\vert S_{1, 2} \vert < r^3$, for sufficiently large $n$ and $p$, with probability tending to $1$.

Reversing the role of $n$ and $p$, and considering $\sum_{j=1}^{(p-1)} \Hdiff{\beta_j}$ instead, one can show that $\vert S_{1,3} \vert < r^3$ for sufficiently large $n$ and $p$ with probability tending to $1$. Thus, combining everything we obtain $\vert S_{1} \vert \leq \vert S_{1,1} \vert + \vert S_{1, 2} \vert + \vert S_{1,3} \vert + \vert S_{1, 4}\vert < 4r^3$ for sufficiently large $n$ and $p$, with probability tending to $1$.

Now, turning our attention to the term $S_2$, we start by writing the expressions for each second order derivative term. Let, $C_\alpha = -(\alpha + 1)(2\pi)^{-\alpha/2}(\sigma^g)^{-(\alpha + 2)}/np$, then
\begin{equation*}
     \E\left[ \Hdifftwo{\lambda} \right]                                                                   
    = C_\alpha \sum_{i=1}^n \sum_{j=1}^p  (u_i^g)^2 (v_j^g)^2 \E\left[ w_{ij} (\alpha Z_{ij}^2 - 1) \right]
    = C_\alpha B_2
\end{equation*}
\noindent where, $B_2 = \E\left[ w_{ij} (\alpha Z_{ij}^2 - 1) \right]$. For the mixed derivative,
\begin{equation*}
      \E\left[ \Hdifftwo{\lambda, \alpha_i} \right] = \sum_{k = 1}^n \E\left[ \Hdifftwo{\lambda, u_k} \right] \partial_{\alpha_i}u_k 
      = \sum_{k = 1}^n C_\alpha \lambda^g \sum_{j=1}^p v_j^g \E\left[ w_{kj}(\sigma^g Z_{kj} + \lambda^g u_k^g v_j^g (\alpha Z_{kj}^2 - 1)) \right] \partial_{\alpha_i}u_k  = 0,
\end{equation*}
\noindent since, $\E(Z_{ij}w_{ij}) = 0$ by symmetry of $g$ and we know $\sum_{i=1}^n u_k^g \partial_{\alpha_i}u_k = 0$. Exchanging the role of $u_i^g$s and $v_j^g$s, we obtain
\begin{equation*}
    \E\left[ \Hdifftwo{\lambda,\beta_j} \right] = 0.
\end{equation*}
\noindent A chain rule of differentiation can be used to obtain the second order derivatives of $\mathcal{H}_{n,p}$ with respect to $\alpha_i$'s as
\begin{equation*}
    \E\left[ \Hdifftwo{\alpha_i, \alpha_j} \right] 
    = \E\left[\sum_{k = 1}^{n} \Hdiff{u_k} \partial^2_{\alpha_i, \alpha_j}u_k + \sum_{k = 1}^n\sum_{l=1}^n \Hdifftwo{u_k, u_l} \partial_{\alpha_i}u_k\partial_{\alpha_j}u_l \right] = \sum_{k=1}^n \E\left[ \Hdifftwo{a_k} \right]  \partial_{\alpha_i}u_k\partial_{\alpha_j}a_k,
\end{equation*}
\noindent since, $\E\left[ \Hdiff{a_k} \right] = 0$ and for $k \neq l$, $\Hdifftwo{a_k, a_l} = 0$. A similar calculation as above reveals that
\begin{equation*}
    \E(\Hdifftwo{a_k})
    = C_\alpha (\lambda^g)^2 B_2.
\end{equation*}
\noindent Combining this with the fact that
\begin{equation*}
    \sum_{k=1}^n \partial_{\alpha_i} u_k \partial_{\alpha_j}a_k
    = \begin{cases}
        (1 - u_n^g)^2 & \text{if, } i = j    \\
        0             & \text{if, } i \neq j
    \end{cases},
\end{equation*}
\noindent yields
\begin{equation*}
    \E\left[ \Hdifftwo{\alpha_i,\alpha_j} \right]
    = \begin{cases}
        C_\alpha (\lambda^g)^2 B_2 (1 - a_n^g)^2 & \text{ if } i = j    \\
        0                                        & \text{ if } i \neq j
    \end{cases}.
\end{equation*}
\noindent Exact same calculation also holds for $\E(\Hdifftwo{\beta_i,\beta_j})$ with $u_n^g$ replaced by $v_p^g$. Because of Assumption~\ref{assum:open-rect}, $(1 - u_n^g)^2$ and $(1 - v_p^g)^2$ can be bounded below by some $\delta > 0$ independent of $n$ and $p$. Also, note that
\begin{multline*}
     \E\left[ \Hdifftwo{\alpha_i, \beta_j} \right]                                                               
    = \sum_{k=1}^n \sum_{l=1}^p \partial_{\alpha_i}u_k \partial_{\beta_j}v_l \E\left[ \Hdifftwo{u_k, v_l} \right] \\
    = \sum_{k=1}^n \sum_{l=1}^p \partial_{\alpha_i}u_k \partial_{\beta_j}v_l C_\alpha \lambda^g  \E\left[ \sigma^g w_{kl}Z_{kl} + \lambda^g u_k^g v_l^g (\alpha Z_{kl}^2 - 1)  \right]
    = 0
\end{multline*}
\noindent which follows from noting that $\E(z_{kl}w^2_{kl}) = 0$ by symmetry of the density function $g$ and $\sum_k u_k \partial_{\alpha_i}u_k = \sum_l v_l \partial_{\beta_j}v_l = 0$. Furthermore, as shown in~\cite{ghosh2013}, the scale and the location estimator become asymptotically uncorrelated for the classical linear regression setup with normally distributed errors. Therefore, we have
\begin{equation*}
    \E\left[ \Hdifftwo{\alpha_i,\sigma^2} \right] = \E\left[ \Hdifftwo{\beta_j,\sigma^2} \right] \\
    = \E\left[ \Hdifftwo{\lambda,\sigma^2} \right] = 0,
\end{equation*}
\noindent and
\begin{equation*}
    \E\left[ \Hdifftwo{\sigma^2} \right] = (2\pi)^{-\alpha/2} (\sigma^g)^{-(\alpha+4)}\\
    \left[
        \dfrac{\alpha(\alpha+2)}{4\sqrt{1+\alpha}} -
        \dfrac{(\alpha + 1)}{2} B_3
        \right] \asymp \sigma^{-(\alpha+4)},
\end{equation*}
\noindent where $B_3 = \E\left( w_{ij} (1 - 2Z_{ij}^2 + \alpha(1 - Z^2_{ij})^2/2)  \right)$. Therefore, if we consider the $(n+p)\times(n+p)$ matrix $\bb{\Psi}_{n,p}$ whose $(k_1, k_2)$-th element is given by $\E(\Hdifftwo{\eta_{k_1},\eta_{k_2}})$, then $\bb{\Psi}_{n,p}$ turns out to be a diagonal matrix with nonzero entries of the order of $(\sigma^g)^{-(\alpha+2)}/np$ and $\sigma^{-(\alpha+4)}$, among which the minimum is at the order of $(\sigma^g)^{-(\alpha+2)}/np$ due to Assumption~\ref{assum:variance-reduce}. Hence, the minimum eigenvalue of $\bb{\Psi}_{n,p}$ is bounded below by $K_1(\sigma^g)^{-(\alpha+2)}/np$ for some positive finite constant $K_1$.

Now, we decompose $S_2$ by considering elements of $\bb{\Psi}_{n,p}$ as follows
\begin{multline*}
    \sum_{k_1, k_2} \Hdifftwo{\eta_{k_1},\eta_{k_2}} (\eta_{k_1} - \eta_{k_1}^g)(\eta_{k_2} - \eta_{k_2}^g)\\
    = \sum_{k_1, k_2} \left[ \Hdifftwo{\eta_{k_1},\eta_{k_2}} - (\bb{\Psi}_{n,p})_{k_1, k_2} \right] (\eta_{k_1} - \eta_{k_1}^g)(\eta_{k_2} - \eta_{k_2}^g) \\
    + \sum_{k_1, k_2} (\bb{\Psi}_{n,p})_{k_1, k_2} (\eta_{k_1} - \eta_{k_1}^g)(\eta_{k_2} - \eta_{k_2}^g)
\end{multline*}
\noindent Here, we can apply an orthogonal transformation on $(\bb{\eta} - \bb{\eta}^g)$ to express it as a linear combination of the eigenvectors of $\Psi_{n,p}$, so that the second term can be made greater than or equal to $K_1(\sigma^g)^{-(\alpha+2)} r^2/np$. Also, it is evident by a generalized version of Khinchin's Law of Large Numbers that the first summation has the expected value equal to $0$. By a similar routine calculation as above, one can show that the variance of the first term also goes to $0$. Therefore, for sufficiently large $n$ and $p$, with probability tending to $1$, $S_2 > (-r^3 + K_1(\sigma^g)^{-(\alpha+2)} r^2/np)$.

Finally, turning to $S_3$, we note that the expected values of the third order derivatives are bounded as shown below.
\begin{align*}
    \left\vert \E\left[ \partial^3_{\lambda}\Hcal_{n,p} \right] \right\vert
     & = M_1 \left\vert \frac{\sigma^{-(\alpha + 3)}}{np} \sum_{i,j} (a_ ib_j)^3 \right\vert,   \\
    \left\vert  \E\left[ \partial^3_{a_i}\Hcal_{n,p} \right] \right\vert
     & = M_2 \left\vert \frac{\sigma^{-(\alpha + 3)}}{np} \sum_{j} (\lambda b_j)^3 \right\vert, \\
    \left\vert \E\left[ \partial^3_{b_j}\Hcal_{n,p} \right] \right\vert
     & = M_3 \left\vert \frac{\sigma^{-(\alpha + 3)}}{np} \sum_{i} (\lambda a_i)^3 \right\vert, \\
    \left\vert \E\left[ \partial^3_{\sigma^2} \Hcal_{n,p} \right] \right\vert
     & = M_4 \sigma^{-(\alpha+6)},
\end{align*}
\noindent where $M_1, M_2, M_3$ and $M_4$ are positive finite constants. The first three among these are $\mathcal{O}(\sigma^{-(\alpha + 2)}/np)$ and the last one is $\mathcal{O}(\sigma^{-(\alpha+6)})$, which follows from Cauchy-Schwartz inequality and the normalization of $u_i$s and $v_j$s. Combining these bounds along with the continuity of the third order derivative of $\Hcal_{n,p}$ and Assumption~\ref{assum:variance-reduce}, we obtain that $\vert S_3\vert \leq M \sigma^{-(\alpha+6)}$ for sufficiently large $n$ and $p$, and for some sufficiently large finite positive constant $M$ independent of $n$ and $p$. Therefore, using the bounds for the individual terms of the Taylor's series, we have
\begin{equation}
    \mathcal{H}_{n,p}(\bb{\eta}) - \mathcal{H}_{n,p}(\bb{\eta}^g) > (-5r^3  \\
    + \frac{K_1}{np} (\sigma^g_{n,p})^{-(\alpha+2)} r^2 - M (\sigma^g_{n,p})^{-(\alpha + 6)} r^3 ),
    \label{eqn:Hcal-lower-bound}
\end{equation}
\noindent with probability tending to $1$ for sufficiently large $n$ and $p$. Now since $(\sigma^g)^4 \asymp (np)^{-1}$ due to Assumption~\ref{assum:variance-reduce} and as $\sigma^g \rightarrow 0$ as $n$ and $p$ tends to infinity, it follows that
\begin{equation*}
    \lim_{n, p \rightarrow \infty} \dfrac{K_1 (\sigma^g_{n,p})^{-(\alpha+2)}/np }{5 + M (\sigma^g_{n,p})^{-(\alpha + 6)}} \\
    = \dfrac{K_1 (\sigma^g_{n,p})^{-(\alpha+6)} }{5 + M (\sigma^g_{n,p})^{-(\alpha + 6)}} = K_2 \in (0, \infty).
\end{equation*}
\noindent Choosing $r < K_2$ ensures that the lower bound in~\eqref{eqn:Hcal-lower-bound} remains positive, i.e., $\mathcal{H}_{n,p}(\bb{\eta}) > \mathcal{H}_{n,p}(\bb{\eta}^g)$ for any $\bb{\eta}$ satisfying $\Vert \bb{\eta} - \bb{\eta}^g \Vert_2 = r$ (where $\Vert \cdot \Vert_2$ denotes the Euclidean $L_2$ norm). This is exactly what we intended to show at the beginning.

Finally, since $\bb{\eta}^\ast$ is consistent for $\bb{\eta}^g$, an application of the continuous mapping theorem completes the proof.

\bibliographystyle{unsrtnat}
\bibliography{reference-theory}  

\end{document}